\documentclass{article}
\usepackage{arxiv}
\usepackage{amsmath}
\usepackage{graphicx}
\usepackage{epstopdf}
\usepackage[utf8]{inputenc} 
\usepackage[T1]{fontenc}    
\usepackage{hyperref}       
\usepackage{url}            
\usepackage{booktabs}       
\usepackage{amsfonts}       
\usepackage{nicefrac}       
\usepackage{microtype}      
\usepackage{lipsum}
\usepackage{algorithm}
\usepackage{algorithmic}
\usepackage{booktabs}
\usepackage{lscape}
\usepackage{caption}
\usepackage{subcaption}
\usepackage{multirow}
\newtheorem{thm}{Theorem}[section]
\newtheorem{lem}[thm]{Lemma}

\newtheorem{defi}{Definition}[section]

\newcommand{\T}{\ensuremath{\mathcal{T}}}

\newcommand{\V}{\ensuremath{\mathcal{V}}}

\newenvironment{Assumptions}
{
\setcounter{enumi}{0}

\begin{enumerate}}
{\end{enumerate} }

\numberwithin{equation}{section} \allowdisplaybreaks

\title{On the Development of a Coupled Non-linear Telegraph-Diffusion Model for Image Restoration}

\author{
  Sudeb Majee \\
  School of Basic Sciences\\ 
  Indian Institute of Technology Mandi\\
  PIN 175005, INDIA\\
  \texttt{sudebmajee@gmail.com} \\
    \And
  Subit K. Jain \\
  School of Basic Sciences\\ 
  Indian Institute of Technology Mandi\\
  PIN 175005, INDIA\\
  \texttt{jain.subit@gmail.com} \\
    \And
  Rajendra K. Ray \\
  School of Basic Sciences\\ 
  Indian Institute of Technology Mandi\\
  PIN 175005, INDIA\\
  \texttt{rajendra@iitmandi.ac.in} \\
    \And
   Ananta K. Majee \\
   Department of Mathematics\\
   Indian Institute of Technology Delhi\\
   PIN 110016, INDIA \\
  \texttt{majee@maths.iitd.ac.in} \\ 
}

\begin{document}
\maketitle

\begin{abstract}
In this work, we propose a telegraph coupled partial differential equation (TCPDE) based model for image restoration. New framework interpolates between a couple of non-linear telegraph equation and a parabolic equation. Proposed strategy can be applied to significantly preserve the oscillatory and texture pattern in an image, even in low signal-to-noise ratio. First, we prove that the present model has a unique global weak solution using Banach's fixed point theorem. Then apply our model over a set of gray-level images to illustrate the superiority of the proposed model over the recently developed hyperbolic-parabolic PDE based models as well as coupled diffusion-based model.
\end{abstract}

\section{Introduction}
\label{intro}
In digital image processing applications, image smoothing is often considered as a significant pre-processing step to make accurate and reliable subsequent image analysis \cite{aubert2006mathematical,scherzer2010handbook}. The principal objective of image denoising algorithms is to achieve the balance between the noise removal and preservation of edges and key features. In the present scenario, partial differential equation (PDE) based approaches are an influential tool for image denoising problem. Due to their well studied mathematical properties and approximation processes, PDE based methods can effectively simulate and preserve slowly varying signals and essential feature of images. Most popular PDE based approaches are anisotropic diffusion models \cite{aubert2006mathematical,jain2015non,jain2016alternative,perona1990scale,weickert1997review,weickert1998anisotropic}, fourth-order PDE based models \cite{deng2019hessian,liu2011adaptive,liu2015fourth,lysaker2003noise,siddig2018image,you1998image,you2000fourth,zhang2017adaptive} total variation models \cite{chambolle1997image,rudin1992nonlinear,thanh2019review,tsai2005total,xu2019novel}, second-order variational models \cite{zanella2018serial,zanetti2016numerical}, and coupled diffusion models \cite{guidotti2015anisotropic,guo2011reaction,acpde2019rgate,nitzberg1992nonlinear}, which are successfully employed to provide a trade-off between edge preservation and noise removal. Recently, Jain et al. \cite{acpde2019rgate} proposed a coupled PDE based diffusion model for additive Gaussian noise removal problem, which takes the following form:
\begin{align}
&I_t = \text{div}(g(u)\nabla I)-2 \lambda v,  \hspace{2.8cm} \text{in}\,\,\, \Omega_T:= \Omega \times (0,T)\,, \label{cpdemaina} \\
&u_t =\kappa\big( h(|\nabla I_{\xi}|^{2})-u+\frac{\nu^2}{2}\Delta u\big),  \hspace{1.7cm} \text{in} \,\,\,\Omega_T\,, \label{cpdemainb} \\
&v_t = \Delta v-(I_0-I),  \hspace{3.6cm} \text{in}\,\,\, \Omega_T \,, \label{cpdemainc} \\
&\partial_n I=0\,,\,\,\, \partial_n u=0\,,\,\,\, v=0, \hspace{2.4cm} \text{on}\,\,\, \partial \Omega_T:=  \partial \Omega \times (0,T)\,, \label{cpdemaind}\\
&I(x,0)=I_0(x)\,,\,\,\, v(x,0)=0\,,\,\,\, u(x,0)=G_{\xi}\ast|\nabla I_0|^{2},\hspace{1.2cm} \text{in} \,\,\,\Omega\,. \label{cpdemaine}
\end{align}
 In the above, $\Omega \subset \mathbb{R}^2$ is the spatial domain of the original image $I$ and the observed noise image $I_0$ and $T>0$ is a prescribed time. Also, $\text{div}$ and $\nabla$ represents the divergence and gradient operator, respectively. $h$ is a smooth version of the truncate function of $|\nabla {I_\xi}|$, where $I_{\xi} =G_{\xi} \ast I,$ ``$\ast$" is the convolution operator, $G_\xi$ is the two dimensional Gaussian kernel. $\Delta$ is the Laplace operator, and $\kappa > 0$, $\nu > 0$ are parameters to be specified and $\lambda$ is the weight parameter calculated as mentioned in \cite{rudin1992nonlinear}. $\partial_n$ denotes the derivative at the boundary surface $ \partial \Omega $ in the outward normal direction n. The diffusion function $g$ is chosen as
\begin{equation}
\label{edge}
g(u)=\dfrac{1}{1+\frac{|u_{\xi}|}{k^2}},
\end{equation}
$k>0$ is the threshold parameter. The equation \eqref{cpdemainb} and \eqref{cpdemainc} were used to achieve a suitable edge map ($u$) and fidelity ($v$) between noisy image and restored image at each scale. Over the last three decades, parabolic PDEs have acquired the center stage in the field of image denoising. Besides the parabolic PDEs, the hyperbolic PDEs which describes oscillations within objects, could also improve the quality of the detected edges more suitable than diffusion based models and so enhance the image better than parabolic PDEs \cite{cao2010class,jain2016edge,averbuch2006edge,ratner2007image,sun2016class,yang2014kernel}. In this regard, Ratner and Zeevi \cite{ratner2007image} proposed a telegraph-diffusion equation (TDE model), which interpolates between the diffusion equation and the wave equation. The TDE model takes the form,
\begin{align*}
&I_{tt}+\gamma I_t  =\text{div}(g(|\nabla I|)\nabla I),  \hspace{0.7cm} \text{in} \hspace{0.2cm} \Omega_T, \\
&\partial_n I=0,      \hspace{3.5cm} \text{in}
\hspace{0.2cm} \partial \Omega_T,\\
&I(x,0)=I_0(x), \hspace{0.2cm} I_t(x,0)=0,	\hspace{0.4cm} \text{in} \hspace{0.2cm}\Omega,
\end{align*}
where  $g(|\nabla I|)=\dfrac{1}{1+\left(\frac{|\nabla I|}{K}\right)^2}$  is an edge-controlled diffusion function which preserves the important features and smoothens the unwanted signals, $K$ is a threshold constant. Here the elasticity and damping
parameters are denoted by $g$ and $\gamma$ respectively. Although the TDE model performs better, it is challenging to confirm the well-posedness of their model. Subsequently, a regularized version of the TDE model has been introduced by Cao et al. \cite{cao2010class}. Their model takes the form 
\begin{align*}
&I_{tt}+\gamma I_t  =\text{div}(g(|\nabla I_{\xi}|)\nabla I),  \hspace{0.5cm} \text{in} \hspace{0.2cm} \Omega_T, \\
&\partial_n I=0,      \hspace{3.4cm} \text{in}
\hspace{0.2cm} \partial \Omega_T,\\
&I(x,0)=I_0(x), \hspace{0.2cm} I_t(x,0)=0,	\hspace{0.4cm} \text{in} \hspace{0.2cm}\Omega.
\end{align*}
The authors replace the gradient $|\nabla I|$  by $|\nabla I_\xi| $ in the edge-controlled function $g$ in the TDE model \cite{ratner2007image}. It has been proved that equations of the form of telegraph equations provide better edge preservation and image enhancement when compared with ordinary diffusion-based methods. For different elasticity coefficients, the telegraph model can be treated as the improved versions of the corresponding nonlinear diffusion models. Although these methods are efficient in the restoration of discontinuous and noisy signals, their performance is not satisfactory in the case of higher noise level or low SNR images.
To overcome this issue, several telegraph models, based on non-linear diffusion method, were proposed \cite{jain2016edge,sun2016class,yang2014kernel}.
In low SNR images, the appropriate separation of noise and important texture information can be viewed as a highly complex problem. The parabolic PDE based restoration methods strongly depend on the diffusion function $g(\cdot)$, to preserve the small variations in the image. In conventional approaches, spatial regularization has been used for diffusion function, which is not able to inject the past information into the diffusion process. To alleviate this shortcoming in the parabolic PDE based models as well as single hyperbolic PDE based models a coupled hyperbolic-parabolic coupled system for image restoration problem was introduced by Sun et al. \cite{sun2016class} which takes the form 
\begin{align}
&I_t-\text{div} \left(g\left(u\right)\nabla I\right)=0,  \hspace{4.0cm} \text{in} \,\,\, \Omega_T, \label{sunmaina} \\
&u_{tt} +u_t-\lambda\text{div} \left(\nabla u\right)-(1-\lambda)\left( \vert \nabla I \vert - u \right) =0,  \hspace{0.6cm} \text{in} \,\,\,\Omega_T\,, \label{sunmainb}  \\
&\partial_n I=0\,,\,\,\, \partial_n u=0\,,\,\,\,  \hspace{4.4cm} \text{on}\,\,\, \partial \Omega_T\,, \label{sunmainc} \\
&I\left(x,0\right)=I_0\left(x\right)\,, \,\,\, u\left(x,0\right)=u_0(x),\,\, u_t\left(x,0\right)=0, \hspace{0.3cm}  \text{in} \,\,\,\Omega\,. \label{sunmaind}
\end{align}
$\lambda>0$ is a balancing parameter.
They considered the following two diffusion functions:
\begin{align*}
g(s)=\dfrac{1}{1+\left( \dfrac{s}{K} \right)^2} \,\,\,\, \text{or} \,\,\,\, g(s)=\vert s \vert^{-1}\,\,\,\, \text{with}\,\,\,\, K>0.
\end{align*}
In the above model, the authors have used telegraph equation only for the edge variable $u$, but it would be superior to use telegraph equation for the image variable $I$ as well as for the edge variable $u$. Also, they have not incorporated the fidelity term \cite{acpde2019rgate} into the model \eqref{sunmaina}-\eqref{sunmaind}, which keeps the restored image close to the original image. 

To overcome these issues, present work aims to systematically develop a new non-linear coupled telegraph diffusion system to deal with the additive Gaussian noise problem. Our motivation is based on the improvement in restoration ability of the proposed coupled system for low SNR images. Inspired by the merit of coupled partial differential equation (ACPDE) model \cite{acpde2019rgate} and telegraph-diffusion equation based models \cite{ratner2007image,sun2016class}, we propose the following telegraph coupled partial differential equation (TCPDE) model:
\begin{align}
&I_{tt} +\alpha I_t = \text{div} \left(g\left(u\right)\nabla I\right)-2 \lambda v,  \hspace{3.0cm} \text{in}\,\,\, \Omega_T\,, \label{maina} \\
&u_{tt} +\beta u_{t} =\kappa\left( h\left(|\nabla I_{\xi}|\right)-u+\frac{\nu^2}{2}\Delta u\right),  \hspace{1.9cm} \text{in} \,\,\,\Omega_T\,, \label{mainb} \\
&v_{t} = \Delta v-\left(I_0-I\right),  \hspace{4.7cm} \text{in}\,\,\, \Omega_T \,, \label{mainc} \\
&\partial_n I=0\,,\,\,\, \partial_n u=0\,,\,\,\, \partial_n v=0, \hspace{3.3cm} \text{on}\,\,\, \partial \Omega_T\,, \label{maine}\\
& \begin{cases}  
I\left(x,0\right)=I_0\left(x\right)\,,\ I_t\left(x,0\right)=0  \,,\,\,\, v\left(x,0\right)=0\,, \\
u\left(x,0\right)=G_{\xi}\ast|\nabla I_0|^{2},\,\, u_t\left(x,0\right)=0,\,\, \hspace{2.1cm} \text{in} \,\,\,\Omega\,. \label{maind}
\end{cases}
\end{align}
where $\alpha, \beta, \kappa, \nu > 0$ are parameters to be specified. In this model, the diffusion coefficient $g(u)$ is chosen same as \eqref{edge}. In this process, the fidelity term between $I$ and $I_0$ can be obtained by function $v$ [see equation \eqref{mainc}], whereas the edge variable $u$ is calculated from equation \eqref{mainb}. The advantages of the proposed model are: (i) parabolic nature of the equation \eqref{maina} remove the noise efficiently, and the hyperbolic nature enhance the image edges better than simple parabolic PDE based models, (ii) hyperbolic nature of the equation \eqref{mainb} detect the image edges better than parabolic PDE,  this extra equation calculates the edge variable $u$ which inject the past information into the diffusion process in equation \eqref{maina}, (iii) fidelity variable $v$ keeps the restored image close to the original image. Overall the proposed approach enables us to provide more flexibility in the diffusion process along the curves of discontinuities. Furthermore, we study the existence and uniqueness of a weak solution of the proposed model using Banach fixed point theorem on an appropriate function space. Moreover, to obtain tangible results, we utilize a robust scheme, which is fast and easy to implement \cite{jain2015comparative,jain2015iterative}. The proposed model has been applied to several natural images. Numerical results illustrate that the proposed algorithm outperforms the existing original Telegraph diffusion model\cite{ratner2007image}, Cao model \cite{cao2010class}, SYS model \cite{sun2016class}, and ACPDE model \cite{acpde2019rgate} in terms of preserving the image structures and noise removal.

The rest of the paper is organized as follows. In section, \ref{sec:analysis}, we study the well-posedness of the proposed model. A numerical realization of the proposed model is shown in section \ref{sec:numerical}. Numerical experiments are carried out and compared with other existing models, in section \ref{sec:results}. Finally, the work is concluded in section \ref{sec:conclusion}.
\section{Well-posedness of weak solutions}
\label{sec:analysis}
In this section we study the well-posedness of the proposed system \eqref{maina}-\eqref{maind} using Banach fixed-point theorem \cite{lcevans1998}. For simplicity we
choose all the constants involved in the equations \eqref{maina}-\eqref{maind} equals to $1$.

\subsection{Technical framework and statement of the main result}
Throughout this section, $C$ denotes a generic positive constant. For $1\le p\le \infty$, we denote by $(L^p, \|\cdot\|_{L^p})$ the standard spaces of $p$-th order integrable
functions on $\Omega$. For $r\in \mathbb{N}$,
we write $(H^r, \|\cdot\|_{H^r})$ for usual Sobolev spaces on $\Omega$, and $H^{-1}=(H^{1}_0)^\prime$. 
We introduce the solution space $W(0,T)= \big(W_1(0,T)\big)^2\times W_2(0,T)$ for the
problem \eqref{maina}-\eqref{maind}, where
\begin{align*}
 W_1(0,T)&=\Big\{w: w\in L^\infty(0,T; H^1)\,, w_t \in L^\infty(0,T; L^2); \,w_{tt} \in L^2(0,T; (H^1)')\Big\}\,, \\
 W_2(0,T)&= \Big\{w: w\in L^\infty(0,T; H^1);\,  w_t \in L^2(0,T; L^2)\Big\}\,.
\end{align*}
Note that the space $W_i(0,T)~(i=1,2)$ is a Hilbert space for the graph norm, see \cite{lions1968controle}. 
\begin{defi}[Weak solution]\label{defi:weak}
A triplet $(I,u,v)$ is called a weak solution of \eqref{maina}-\eqref{maind} if
\begin{itemize}
 \item[a)] $I,u \in W_1(0,T),~v\in W_2(0,T)$ and \eqref{maind} holds.
 \item[b)] For all $\phi \in H^1$, $\psi \in H_0^1$, and a.e $t\in (0,T)$, there hold
 \begin{align*}
&  \big\langle  I_{tt}, \phi \big\rangle + \int_{\Omega} I_t \phi \,dx + \int_{\Omega} g(u)\nabla I\cdot \nabla \phi\,dx = 2 \int_{\Omega} v \phi\,dx\,, \\
&   \big\langle  u_{tt}, \phi\big\rangle + \int_{\Omega} u_t \phi \,dx + \frac{1}{2}\int_{\Omega} \nabla u\cdot \nabla \phi\,dx 
   + \int_{\Omega}u\phi\,dx = 2 \int_{\Omega} h(|\nabla I_\xi|) \phi\,dx\,, \\
&   \int_{\Omega} v_t \psi \,dx + \int_{\Omega} \nabla v\cdot \nabla \phi\,dx = - \int_{\Omega} (I_0-I) \psi\,dx\,.
 \end{align*}
\end{itemize}
\end{defi}
 As we mentioned, our aim is to establish well-posedness of weak solutions of \eqref{maina}-\eqref{maind}, and we will do so under the following assumptions:
 \begin{Assumptions}
  \item \label{A1} $h:\mathbb{R}^+ \rightarrow \mathbb{R}^+$ is a bounded, Lipschitz continuous function with Lipschitz constant $c_h$ such that
\begin{align*}
 0 \leq h(\tilde{u})\leq 1\, \quad \forall\, \tilde{u}\in \mathbb{R}^+\,.
 \end{align*}
Moreover, $h^\prime$ is Lipschitz continuous with Lipschitz constant $c_{h^\prime}$.
\item \label{A2} $I_0 \in H^2$.
 \end{Assumptions}
We observe from \eqref{edge} that, $g: \mathbb{R} \rightarrow \mathbb{R}^+ $ is a bounded, decreasing and Lipschitz continuous function with Lipschitz constant $\frac{C_{\xi}}{k^2}$.
Moreover, $g(0)=1 $ and $ \underset{u \rightarrow + \infty }\lim g(u)=0 $. 
We are now ready to state the main results of this paper.

\begin{thm}\label{thm:weak}
Under the assumptions \ref{A1}-\ref{A2}, the problem \eqref{maina}-\eqref{maind} admits a unique weak solution $(I,u,v)$ in the sense of Definition \ref{defi:weak}.
\end{thm}
Following \cite{zheng1995nonlinear}, we first prove a local well-posedness result, and then establish a uniform a-priori estimate for the solution. Using uniform moment estimates and continuation method, we 
prove a global existence result.

\subsection{Linearized problem and its well-posedness}
For any positive constants $M_1$, $M_2$, and $M_3$, we define the following convex set $\mathcal{B}^{\T}$ with fixed $\T>0$:
\begin{align*}
\mathcal{B}^{\T}=
\begin{cases}
\left(u,u_t,u_{tt}\right) \in \mathcal{B}^{\T}_1 \equiv L^{\infty}\left(0,\T;H^1\right) \times L^{\infty}\left(0,\T;L^2\right) \times L^{\infty}\left(0,\T;L^2\right) , \\
\left(v,v_t\right) \in \mathcal{B}^{\T}_2 \equiv L^{\infty}\left(0,\T;H^1\right) \times L^{2}\left(0,\T;L^2\right), \\
\left(I,I_t\right) \in \mathcal{B}^{\T}_3 \equiv L^{\infty}\left(0,\T;H^1\right) \times L^{\infty}\left(0,\T;L^2\right),\\
\Vert u \Vert_{L^{\infty}\left(0,\T;L^2\right)} + \Vert u_t \Vert_{L^{\infty}\left(0,\T;L^2\right)} + \Vert u_{tt} \Vert_{L^{\infty}\left(0,\T;L^2\right)}\leq M_1,\\
\Vert v \Vert_{L^{\infty}\left(0,\T;L^2\right)} +  \displaystyle \int_{0}^{\T}\Vert v_t \Vert^2 dt \leq M_2,\\
\Vert I \Vert_{L^{\infty}\left(0,\T;L^2\right)} +  \Vert I_t \Vert_{L^{\infty}\left(0,\T;L^2\right)} \leq M_3.\\
\end{cases}
\end{align*}

 For any fixed $\left(\bar{I} ,\bar{u},\bar{v}\right) \in \mathcal{B}^{\T} $, consider the following linearized problem:
\begin{align}
&I_{tt} + I_t - {\rm div} \left(g\left(\bar{u}\right)\nabla I\right)=-2 \bar{v}, \hspace{2cm} \text{in}\,\,\,\,\, \Omega_T \label{linmaina} \\
& u_{tt} + u_t-\frac{1}{2}\Delta u + u =    h(|\nabla \bar{I}_{\xi}|),\hspace{1.8cm}\text{in} \,\,\,\Omega_T,\, \label{linmainb} \\
& v_t -\Delta v =\bar{I}-I_0, \hspace{4.3cm} \text{in}\,\,\, \Omega_T \,, \label{linmainc}
\end{align}
with the initial conditions \eqref{maind}. Since $I_0\in H^1$, by using the properties of convolution ($L_2$-estimate), one can easily check that the followings hold, see \cite{cao2010class}
\begin{align}
\gamma: = \dfrac{1}{1+\frac{C_{\xi}\sqrt{M_1}}{k^2}} \le & \bar{g}\le 1 \,, \quad  \left\vert  \bar{g}_t \right\vert \le  \dfrac{C_{\xi}}{k^2}M_1 \,, \label{bound:bar-g-h-derivatives}
\end{align}
where $\bar{g}=g(\bar{u})$. 
Hence by classical Galerkin approximation, there exists a unique solution $\left(I,u,v \right) \in \mathcal{B}^{\T}_3 \times \mathcal{B}^{\T}_1 \times \mathcal{B}^{\T}_2 $
of the linearized problem \eqref{linmaina}-\eqref{linmainc} with the 
initial conditions \eqref{maind}. Moreover, $(I,u,v)$ satisfies the following estimates.

\begin{lem}\label{lem:apriori-estimate-linear}
 The unique solution $(I,u,v)\in \mathcal{B}^{\T}$ of the linearized problem \eqref{linmaina}-\eqref{linmainc} with the 
initial conditions \eqref{maind} satisfies the following: there exists a constant $C>0$, depending only on $\T, \xi, M_1, M_2, M_3, |\Omega|$ and $\|I_0\|$, such that
\begin{itemize}
 \item [i)] $\|I\|_{L^\infty(0,\T; H^2)}^2 + \|I_t\|_{L^\infty(0,\T; L^2)}^2 \le C$,
 \item[ii)] $\|u\|_{L^\infty(0,\T; H^1)}^2 + \|u_t\|_{L^\infty(0,\T; L^2)}^2 + \|u_{tt}\|_{L^\infty(0,\T; L^2)}^2\le C$,
 \item[iii)] $\|v\|_{L^\infty(0,\T; H^1)}^2 + \|v_t\|_{L^2(0,\T; L^2)}^2 \le C$.
\end{itemize}
\end{lem}
\textbf{Proof:}
\noindent{Proof of ${\rm i)}:$}
Multiplying \eqref{linmaina} by $ I_t $ and integrating by parts over $\Omega$ and using Cauchy-Schwarz and  Young's inequalities, we obtain
\begin{align}\label{eq:step 1 I}
&\frac{1}{2}\frac{d}{dt}\|I_t\|_{L^2}^2 +\|I_t\|_{L^2}^2 + \int_{\Omega} \bar{g}\,\nabla I\cdot \nabla I_t\,dx 
\leq \|\bar{v}\|_{L^2}^2 + \|I_t^2\|_{L^2}^2\,.
\end{align}
Note that, thanks to \eqref{bound:bar-g-h-derivatives}
\begin{align}\label{eq:g_delI_delIt}
\int_{\Omega} \bar{g}\, \nabla I\cdot \nabla I_t \,dx
\ge \frac{1}{2}\dfrac{d}{dt}\int_{\Omega} \bar{g}|\nabla I|^2 dx- \frac{C_{\xi M_1}}{2k^2}\|\nabla I\|_{L^2}^2\,.
\end{align}
Combining \eqref{eq:step 1 I} and \eqref{eq:g_delI_delIt}, along with \eqref{bound:bar-g-h-derivatives},  we have
\begin{align}
\dfrac{d}{dt}\|I_t\|_{L^2}^2 + \dfrac{d}{dt}\int_{\Omega} \bar{g}|\nabla I|^2\, dx 
 \le C\Big(\|I_t\|_{L^2}^2 + \int_{\Omega}\bar{g}|\nabla I|^2\,dx \Big)+ 2\|\bar{v}\|_{L^2}^2\,. \label{esti:I-require}
\end{align}
Using Gronwall's inequality, we obtain, for a.e. $t\in (0,\T)$,
\begin{align*}
& \left\Vert I_t \right\Vert^2_{L^2}+ \int_{\Omega}\bar{g}|\nabla I|^2 dx 
 \le  e^{Ct} \left(C_1 + t\,C_2(M_2) \right)\,, \notag \\
& \left\Vert \nabla I \right\Vert^2_{L^2}\le \frac{1}{\gamma}\int_{\Omega}\bar{g}|\nabla I|^2\,dx
\le \frac{1}{\gamma}\ e^{Ct} \left(C_1 + t\,C_2(M_2) \right)\,.
\end{align*}
Thus, one has, for a.e. $t\in (0,\T)$
\begin{align}\label{bound_I_t_nabla_I}
\left\Vert I_t \right\Vert^2_{L^2} + \left\Vert \nabla I \right\Vert^2_{L^2} \leq \widetilde{M_1}e^{Ct} \left(C_1 + t\,C_2(M_2) \right)\,,
\end{align}
where $\widetilde{M_1}= \max\{\gamma^{-1},1 \}$. By using the identity $I(t,x)=I_0(x)+\int_{0}^{t} I_t(s,x)\,ds$, Young's inequality and \eqref{bound_I_t_nabla_I}, 
 we get $\left\Vert I(t) \right\Vert^2_{L^2} \leq \big(2\left\Vert I_0 \right\Vert^2_{L^2} + tC^{\prime}_1\big)e^{Ct} $ which then implies for a.e. $t\in (0,\T)$
\begin{align}\label{boundIH1boundItL2}
\left\Vert I_t \right\Vert^2_{L^2} + \left\Vert I(t) \right\Vert^2_{H^1} \leq e^{Ct} \left(C^{\prime}_2 + tC^{\prime}_3 \right)\,,
\end{align}
where, $C^{\prime}_2=\widetilde{M_1}C_1+2\left\Vert I_0 \right\Vert^2_{L^2}$ and $C^{\prime}_3=\widetilde{M_1} \big(C_2(M_2)+C^{\prime}_1\big)$ with 
$C^{\prime}_1= \frac{2\widetilde{M_1}}{C}\big( C_1 + \T C_2(M_2)$. 
\vspace{.1cm}

We now show that $I\in L^\infty(0,\T; H^2)$, which will play an essential role in the later analysis.  To do so, we follow the arguments as in \cite{lions2012non}. 
Differentiate the equation \eqref{linmaina} w.r.t time, multiply the resulting equation by $I_{tt},$  integrate over $\Omega$, use the inequality 
$\int_{\Omega} \bar{g}\nabla I_t.\nabla I_{tt} dx \geq  \frac{1}{2}\dfrac{d}{dt}\int_{\Omega} \bar{g}|\nabla I_t|^2 dx- \frac{C_{\xi M_1}}{2k^2}\int_{\Omega}|\nabla I_t|^2 dx$ and then
integrate w.r.t time variable from $0$ to $t$ which yields 
\begin{align*}
 & \left\Vert I_{tt} \right\Vert^2_{L^2} + \int_{\Omega} \bar{g}|\nabla I_t|^2 \,dx \notag \\
 & \leq  3\frac{C_{\xi M_1}}{k^2}\int_{0}^{t}\int_{\Omega}|\nabla I_t|^2 \,dx\,dt -2\int_{\Omega} \bar{g}_t \nabla I\cdot \nabla I_{t}\,dx  
+2\int_{0}^{t} \int_{\Omega} \bar{g}_{tt} \nabla I\cdot \nabla I_{t}\,dx\,dt-2\int_{0}^{t}\int_{\Omega}  \bar{v}_t I_{tt}\, dx\,dt \\
& \leq 3\frac{C_{\xi M_1}}{\gamma k^2}\int_{0}^{t} \Big(\int_{\Omega}\bar{g}|\nabla I_t|^2\,dx \Big)\,ds +\dfrac{1}{\epsilon}\|\bar{g}_t\|_{L^\infty} \|\nabla I\|^2_{L^2}
+ \epsilon\|\bar{g}_t\|_{L^\infty}\|\nabla I_{t} \|^2_{L^2}  \nonumber \\
&\hspace{1cm}+ \|\bar{g}_{tt}\|_{L^\infty}\int_{0}^{t} \big(\|\nabla I \|^2_{L^2} + \|\nabla I_{t} \|^2_{L^2}\big)\,ds
+\int_{0}^{t}  \|\bar{v}_t\|^2_{L^2}\,ds + \int_{0}^{t} \|I_{tt}\|^2_{L^2}\,ds \\
& \le \dfrac{1}{\epsilon} \|\bar{g}_t\|_{L^\infty} \|\nabla I\|^2_{L^2}+\|\bar{g}_{tt}\|_{L^\infty} \int_{0}^{t}\|\nabla I \|^2_{L^2} \,ds + 
\int_{0}^{t}  \|\bar{v}_t\|^2_{L^2}\,ds \\
& \hspace{1cm} +\Big( 3\frac{C_{\xi M_1}}{\gamma k^2} + \dfrac{\|\bar{g}_{tt}\|_{L^\infty}}{\gamma}\Big) \int_{0}^{t} \Big(\int_{\Omega}\bar{g}|\nabla I_t|^2\, dx \Big)\,ds 
 + \int_{0}^{t} \|I_{tt}\|^2_{L^2}\,ds + \dfrac{\epsilon \|\bar{g}_t\|_{L^\infty}}{\gamma}\int_{\Omega} \bar{g}|\nabla I_t|^2\,dx \,.
\end{align*}
Observe that $\|\bar{g}_t\|_{L^\infty}, \|\bar{g}_{tt}\|_{L^\infty} \le C(\xi, M_1, \|I_0\|)$. 
Now by the proper choice of $\epsilon$, we can rewrite the above inequality as
\begin{align*}
X(t) \leq Y(t)+ C\int_{0}^{t}X(s)\,ds\,, 
\end{align*}
where
\begin{align*}
\begin{cases}
X(t) = \|I_{tt}\|^2_{L^2} + \int_{\Omega} \bar{g}|\nabla I_t|^2\,dx\,, \notag \\
Y(t) = \dfrac{1}{\epsilon} \|\bar{g}_t\|_{L^\infty} \|\nabla I\|^2_{L^2}+\|\bar{g}_{tt}\|_{L^\infty} \int_{0}^{t}\|\nabla I \|^2_{L^2} \,ds + 
\int_{0}^{t}  \|\bar{v}_t\|^2_{L^2}\,ds\,.
\end{cases}
\end{align*}
An application of Gronwall's lemma together with \eqref{boundIH1boundItL2} yields
\begin{align*}
 \|I_{tt}\|_{L^{\infty}(0,\T;L^2)}+ \|\nabla I_{t}\|_{L^{\infty}(0,\T;L^2)} \leq C\,.
\end{align*}
Since $\nabla \bar{g} \in L^{\infty}\left(0,\T;L^{\infty} \right)$, from \eqref{linmaina}, it is easily to show that  $I \in L^{\infty}(0,\T;H^{2})$, and hence ${\rm i)}$ holds true. 
\vspace{.1cm}

 \noindent{Proof of ${\rm ii)}:$} We multiply \eqref{linmainb} by $u_t $, integrate by parts over $\Omega$ and use Cauchy-Schwarz and
 Young's inequalities. The result is
\begin{align*}
\frac{1}{2}\dfrac{d}{dt}\|u_t\|_{L^2}^2 + \|u_t\|_{L^2}^2 + \frac{1}{4}\dfrac{d}{dt}\|\nabla u\|_{L^2}^2+ \frac{1}{2}\dfrac{d}{dt}\|u\|_{L^2}^2 \leq &\frac{1}{2}\|\bar{h}\|_{L^2}^2 +\frac{1}{2} \|u_t\|_{L^2}^2 \le 
2\vert \Omega \vert + \frac{1}{2} \|u_t\|_{L^2}^2\,.
\end{align*}
Integrating between the time interval form $0$ to $t$, we obtain, for a.e. $t\in (0,\T)$
\begin{align*}
\Vert u(t) \Vert^2_{H^1}+ \Vert u_t \Vert^2_{L^2} \le C_3+2t\,\vert \Omega \vert \,.
\end{align*}
Note that $\bar{h}_t \in L^{\infty}\left(0,\T;L^2\right)$, and hence by regularity theory \cite{lcevans1998}, one can easily prove that $u_{tt} \in L^{\infty}\left(0,\T;L^2\right),$ with
the estimate
\begin{align*}
\|u_{tt}\|_ {L^{\infty}\left(0,\T;L^2\right)} \leq C\big( \|I_0\|, M_1, M_2, \xi \big)\,.
\end{align*}
This shows that ${\rm ii)}$ holds as well. 
\vspace{.1cm}

\noindent{Proof of ${\rm iii)}:$} Multiplying \eqref{linmainc} by $v_t $, integrating by parts over $\Omega$, using Cauchy-Schwarz and Young's inequalities, and then integrating between the time interval
form $0$ to $t$ of the resulting inequality, we obtain, for a.e. $t\in (0,\T)$
\begin{align*}
\Vert \nabla v \Vert^2_{L^2}+\displaystyle\int_{0}^{t}\Vert v_t \Vert^2_{L^2} ds \le  C_4+ t\,M_3\,.
\end{align*}
Again, multiplying \eqref{linmainc} by $v$ and integrating over $\Omega$ and using Cauchy-Schwarz and Young's inequalities along with Gronwall's lemma,  we get
$\Vert v(t) \Vert^2_{L^2}
\le  e^{t}\big(\widetilde{C_4}+tM_3\big)$, and hence
\begin{align}\label{bound-vh1-vt}
\Vert v(t) \Vert^2_{H^1}+\displaystyle\int_{0}^{t}\Vert v_t \Vert^2_{L^2} ds \le e^{t}\left( C_5+ tM_3\right) ,\,\,\, \text{a.e.} \,\, t\in (0,\T).
\end{align}
This finishes the proof. 

\subsection{Well-posedness of a local solution} 
We show that the problem  \eqref{maina}-\eqref{maind} admits a unique solution on a small  time interval. Indeed, we have the following lemma. 
\begin{thm}\label{thm:wellposedness-local-solution}
 There exists a positive time $\T\in (0,T]$, depending only on the data $I_0, h$ and $G_{\xi}$, such that the problem \eqref{maina}-\eqref{maind} admits a unique solution 
 $(I,u,v)$ in $\Omega_{\T}$. Moreover, we have 
 \begin{align*}
  \begin{cases}
 I\in L^\infty(0,\T; H^2),~ I_t \in L^\infty(0,\T;L^2)\,, \\
u\in L^\infty(0,\T; H^1),~  u_t\in L^\infty(0,\T; L^2),~  u_{tt}\in L^\infty(0,\T; L^2)\,, \\
 v\in L^\infty(0,\T; H^1),~  v_t\in L^2(0,\T; L^2)\,.
  \end{cases}
 \end{align*}
\end{thm}
\textbf{Proof:}
In view of Lemma \ref{lem:apriori-estimate-linear}, we see that for small $t$ and hence for small $\T$, the solution 
$\left(I,u,v \right) \in \mathcal{B}^{\T},$ and hence the mapping $\left(\bar{I},\bar{u},\bar{v} \right) \mapsto \left(I,u,v \right)$ maps $\mathcal{B}^{\T}$ into itself. Well-posedness of
the solution of \eqref{maina}-
\eqref{maind} on the time interval $[0,\T]$
would then follows from the Banach fixed point theorem once we establish that the mapping $\left(\bar{I},\bar{u},\bar{v} \right) \mapsto \left(I,u,v \right) $ is a contraction.
\vspace{.1cm}

For fixed $\big(\bar{I},\bar{u},\bar{v} \big),~\big(\bar{\bar I},\bar{\bar u},\bar{\bar v} \big)\in \mathcal{B}^{\T}$, let $(I,u,v)$ and $\big(\tilde{I},\tilde{u},\tilde{v}\big)$ be the corresponding solutions of the linearized problem \eqref{linmaina}-\eqref{linmainc}.
Let us denote $\left(\zeta,\theta,\eta\right)=\big(I-\tilde{I},u-\tilde{u},v-\tilde{v}  \big)$ and $
\left(\bar{\zeta},\bar{\theta},\bar{\eta}\right)=\big(\bar{I}-\bar{\bar I},\bar{u}-\bar{\bar u},\bar{v}-\bar{\bar v} \big)$. 
Consider the equation for $\theta,$ i.e.
\begin{align}\label{thetaequation}
\theta_{tt}+\theta_t - \frac{1}{2}\Delta \theta + \theta = \big( \bar h - \bar{\bar h} \big)\,,
\end{align}
where $\bar h= h\big(|\nabla \bar{I}_{\xi}|\big) $ and $\bar{\bar h}=h\big(|\nabla \bar{\bar{I}}_{\xi}|\big).$ Like an analogous way to the estimates established in Lemma \ref{lem:apriori-estimate-linear}, we 
multiply \eqref{thetaequation} by $ \theta_{t}$, integrate by parts over $\Omega$, and then use Young's inequality to obtain
\begin{align*}
 \|\theta_t\|_{L^2}^2+ \frac{d}{dt}\Big( \|\theta_t\|_{L^2}^2+ \frac{1}{2} \|\nabla \theta\|_{L^2}^2 + \|\theta\|_{L^2}^2\Big) \leq &\|\bar h - \bar{\bar h}\|_{L^2}^2\,.
\end{align*}
Thanks to Lipschitz continuity of $h$ and Young's inequality for convolution, we see that \\
$\|\bar h - \bar{\bar h}\|_{L^2}^2 \le C(c_h, \xi)\|\bar{\zeta}\|_{L^2}^2$, 
and hence 
\begin{align*}
 \|\theta_t\|_{L^2}^2 +\|\theta\|_{H^1}^2 \le C\,t\,\underset{0 \leq t \leq \T}{\text{sup}}  \|\bar{\zeta}\|_{L^2}^2.
\end{align*}
We would like to estimate $\|\theta_{tt}\|_{L^2}$. We use standard methodology \cite{lcevans1998}, i.e., differentiate \eqref{thetaequation} w.r.t time, multiply the resulting equation by $\theta_{tt}$,
integrate over $\Omega$, and then use
Young's inequality. The result is
\begin{align*}
&\dfrac{d}{dt}\left(\Vert \theta_{tt} \Vert^2_{L^2} +\Vert \theta_t \Vert^2_{H^1} \right) \leq  \|\bar{h}_t-\bar{\bar h}_t \|^2_{L^2}.
\end{align*}
In view of Lipschitz continuity of $h^\prime$, one has  $\|\bar{h}_t-\bar{\bar h}_t \|^2_{L^2} \le C \big( \|\bar{\zeta}\|^2_{L^2}  +  \|\bar{\zeta}_t  \|^2_{L^2} \big)$, and hence
\begin{align*}
&\|\theta_{tt} \|^2_{L^2}  \leq C\,t\,\sup_{0 \leq t \leq \T}\big( \|\bar{\zeta}\|^2_{L^2}  +  \|\bar{\zeta}_t  \|^2_{L^2} \big)\,.
\end{align*}
Thus, we have
\begin{align}\label{reg_theta-final}
&\|\theta_t\|_{L^2}^2 +\|\theta\|_{H^1}^2 + \|\theta_{tt} \|^2_{L^2}  \leq C\,t\,\sup_{0 \leq t \leq \T}\big( \|\bar{\zeta}\|^2_{L^2}  +  \|\bar{\zeta}_t  \|^2_{L^2} \big)\,.
\end{align}
To derive estimate for the solution $\eta$, a weak solution of the PDE $~~
\eta_t - \Delta \eta = \big( \bar I - \bar{\bar I}\big)$, 
we proceed similarly as in the derivation of \eqref{bound-vh1-vt}, and obtain 
\begin{align}\label{reg_eta-final}
\|\eta(t) \|^2_{H^1}+\int_{0}^{t}\|\eta_t \|^2_{L^2} \,ds \le  t\,e^{t}\,\sup_{0 \leq t \leq \T}\|\bar{\zeta}\|^2_{L^2}\,.
\end{align}
Now focus on the equation for $\zeta,$ i.e.
\begin{align}\label{zetaequation}
\zeta_{tt}+\zeta_{t}-\nabla\left( \bar{g} \nabla \zeta \right)=\nabla\big(( \bar{g}-\bar{\bar g})\nabla\tilde{I} \big)-2\bar{\eta}\,,
\end{align} 
where, $\bar g= g(\bar {u}) $ and $\bar{\bar g}=g(\bar{\bar u})$. We multiply \eqref{zetaequation} by $ \zeta_{t} $ and integrate
over $\Omega$ to have 
\begin{align}
\frac{d}{dt} \big( \|\zeta_{t} \|_{L^2}^2  + \int_{\Omega} \bar{g}|\nabla \zeta|^2\,dx \big)
& \leq 2||\zeta_{t}||^2_{L^2} + C \|\nabla \zeta\|_{L^2}^2 + \widetilde{C}_1 || \bar{g}-\bar{\bar g}||^2_{L^\infty}+ \widetilde{C}_2 ||\nabla \left(\bar{g}-\bar{\bar g}\right)||^2_{L^\infty}+||\bar{\eta}||^2_{L^2}\,, \label{step 1 zetaequation}
\end{align}
where $\widetilde{C}_1=||\Delta \tilde{I}||^2_{L^{\infty}(0,\T; L^2)} $ and $\widetilde{C}_2=||\nabla \tilde{I}||^2_{L^{\infty}(0,\T; L^2)}.$ By Lipschitz continuity of $g$,
$L^p$-estimate for convolution and \eqref{bound:bar-g-h-derivatives}
along with H\"{o}lder's inequality, we see that 
\begin{equation*}
||\bar{g}-\bar{\bar g}||_{L^\infty} \leq C||\bar{u}-\bar{\bar u}||_{L^2}; \quad ||\nabla \left(\bar{g}-\bar{\bar g}\right)||^2_{L^\infty} \le
\widetilde{C}_3 \Vert \bar{\theta} \Vert^2_{L^2}\,,
\end{equation*}
where $ \widetilde{C}_3 = \Big( \dfrac{2}{k^2} \|\nabla G_{\xi} \|^2_{L^2} + \dfrac{2}{k^2} \|({\bar g} + \bar{\bar g}) \|^2_{L^\infty}
\|\nabla G_{\xi} \ast\bar{\bar u} \|^2_{L^\infty} \Big).$
Hence from \eqref{step 1 zetaequation}, we have
\begin{align*}
\|\zeta_{t} \|_{L^2}^2 + \|\nabla \zeta\|_{L^2}^2 
\leq C\, t\,e^{Ct}\sup_{0 \leq t \leq \T} \big(\|\bar{\theta}\|^2_{H^1}+\|\bar{\eta}\|^2_{H^1}\big)\,. 
\end{align*}
 Moreover, one can easily show that
\begin{align}\label{reg_zeta-final}
\|\zeta_{t} \|_{L^2}^2 + \|\zeta \|_{H^1}^2 
&\leq  C t\,e^{Ct} \sup_{0\le t\le \T}\big(||\bar{\theta}||^2_{H^1}+\|\bar{\eta}\|^2_{H^1}\big)\,. 
\end{align}

Combining \eqref{reg_theta-final}, \eqref{reg_eta-final} and \eqref{reg_zeta-final}, we have
\begin{align}\label{contraction_small_t}
&\Vert \zeta_{t} \Vert_{L^2}^2 + \Vert \zeta \Vert_{H^1}^2+\Vert \theta_t \Vert^2_{L^2}+\Vert \theta \Vert^2_{H^1} + \|\theta_{tt}\|_{L^2}^2+\Vert \eta \Vert^2_{H^1}+\displaystyle\int_{0}^{t}\Vert \eta_t \Vert^2_{L^2} ds \nonumber \\
&\leq C\,t\,e^{Ct}\sup_{0 \leq t \leq \T}\Big\{\|\bar{\zeta}_t\|^2_{L^2}+\|\bar{\zeta} \|_{H^1}^2 + \|\bar{\theta}\|^2_{H^1}+\|\bar{\theta}_t \|^2_{L^2}+ \|\bar{\theta}_{tt} \|^2_{L^2} +\|\bar{\eta}\|^2_{H^1}
+\int_{0}^{t}\|\bar{\eta}_t\|^2_{L^2}\,ds \Big\}\,.
\end{align}
Hence the contraction property of the mapping $\left(\bar{I},\bar{u},\bar{v} \right) \mapsto \left(I,u,v \right)$, for $t$ small enough, follows immediately from \eqref{contraction_small_t}. This shows that,
there exists a unique solution $(I,u,v)$ of the underlying problem \eqref{maina}-\eqref{maind} over a time interval $(0,\T)$ for small $\T>0$. This completes the proof.
\begin{lem}[Uniform a priori estimate]\label{lem:uniform-apriori-estimate}
There exists a positive constant $C_T=C(I_0, G_{\xi}, T, h)$ such that the solution $(I,u,v)$ of the underlying problem \eqref{maina}-\eqref{maind} verifies the following estimate:
 \begin{align}\label{uniform_final}
&\sup_{\tau \in (0,T]} \Big\{\Vert I_t(\tau) \Vert^2_{L^2} +\Vert I(\tau) \Vert^2_{H^1}+\Vert u_t(\tau) \Vert^2_{L^2} + \Vert u(\tau) \Vert^2_{H^1} \notag \\
& \hspace{2cm} +\|u_{tt}(\tau)\|^2_{L^2} +\|v(\tau)\|^2_{H^1} +\int_{0}^{\tau} \|v_t(s)\|^2_{L^2}\,ds \Big\} \leq C_{T}\,.
\end{align}
\end{lem}

\textbf{Proof:}
First consider the equation \eqref{maina}. Arguing similarly as in the derivation of \eqref{esti:I-require}, we get
\begin{align*}
&\dfrac{d}{dt}\Big(\Vert I_t \Vert^2_{L^2} + \int_{\Omega} g|\nabla I|^2\, dx \Big)\leq C\Big(\Vert I_t \Vert^2_{L^2} + \int_{\Omega}g|\nabla I|^2\,dx\Big)+ 2\Vert v \Vert^2_{L^2}\,.
\end{align*}
 Since $
\dfrac{d}{dt}\|I\|_{L^2}^2 = \int_{\Omega}2I I_t \,dx \le \|I\|^2_{L^2} + \|I_t \|^2_{L^2} $, 
we obtain
\begin{align}
 &\dfrac{d}{dt}\Big(\|I_t\|^2_{L^2} + \int_{\Omega} g|\nabla I|^2\, dx  + \|I\|_{L^2}^2\Big)\leq C \Big(\|I_t\|^2_{L^2} + \int_{\Omega}g|\nabla I|^2\,dx\Big)+ 2\|v\|^2_{L^2}
 + \|I\|^2_{L^2}\,. \label{eq:uniform-I}
\end{align}
Multiplying \eqref{mainb} by $u_t$, and applying Cauchy-Schwarz and Young's inequalities, we have
\begin{align}\label{eq:uniform-u}
\dfrac{d}{dt}\left(2\,\Vert u_t \Vert^2_{L^2}+\Vert \nabla u \Vert^2_{L^2} + 2\,\Vert u \Vert^2_{L^2} \right) \leq & 2\vert \Omega \vert.
\end{align}
Moreover, upon differentiating \eqref{mainb} w.r.t time and then tested with $u_{tt}$, one has
\begin{align}
 \dfrac{d}{dt} \big( \|u_{tt}\|^2_{L^2} + \|u_t\|^2_{H^1}\big) \le \|h_t\|_{L^2}^2 \le C|\Omega|\, .\label{eq:unifor_u_tt_ext}
\end{align}
Next we multiply \eqref{mainc} by $v_t$ and $v$ resp. and then use integration by parts over $\Omega$ along with Cauchy-Schwarz and Young's inequalities to arrive at 
\begin{equation}\label{eq:uniform-v}
\begin{aligned}
&\Vert v_t \Vert^2_{L^2} + \dfrac{d}{dt} \Vert \nabla v \Vert^2_{L^2} \leq  2\,\Big(\Vert I_0\Vert^2_{L^2} +\Vert I \Vert^2_{L^2}\Big)\,,\\
&\dfrac{d}{dt} \Vert v \Vert^2_{L^2} \leq \Vert v \Vert^2_{L^2}+ \Vert I_0\Vert^2_{L^2} +\Vert I \Vert^2_{L^2}\,.
\end{aligned}
\end{equation}
Combining \eqref{eq:uniform-I}, \eqref{eq:uniform-u}, \eqref{eq:unifor_u_tt_ext}, and \eqref{eq:uniform-v}, we have 
\begin{align}
\dfrac{d}{dt}p(t) \leq - \|v_t\|^2_{L^2} + (1+C)|\Omega| +3\|I_0\|^2_{L^2} + C p(t)\,, \label{uni-1}
\end{align}
where
\begin{align*}
p(t)=\|I_t\|^2_{L^2} + \int_{\Omega} g|\nabla I|^2\,dx + \|I\|^2_{L^2}+ \|u\|^2_{H^1} + 2\|u_t\|^2_{L^2} + \|u_{tt}\|^2_{L^2} + \|v\|^2_{H^1}.
\end{align*}
From \eqref{uni-1}, we get
\begin{align*}
\dfrac{d}{dt}\left\{e^{-Ct}p(t)\right\}+e^{-Ct}\|v_t\|^2_{L^2} \leq C\big(|\Omega | + \|I_0\|^2_{L^2}\big)\,.
\end{align*}
Integration for $ t\in (0,\tau),$ for any $\tau \in (0,T]$, yields the following inequality at time $\tau:$
\begin{align*}
&\|I_t\|^2_{L^2} + \int_{\Omega} g|\nabla I|^2 \,dx+\|I\|^2_{L^2}+\|u\|^2_{H^1} + 2\|u_t\|^2_{L^2} + \|u_{tt}\|^2_{L^2}
+\|v\|^2_{H^1}+\int_{0}^{\tau}\|v_t\|^2_{L^2}\, ds \nonumber\\
 & \hspace{5cm}\quad \leq e^{C\tau}\Big\{p(0)+ C\,\int_{0}^{\tau}\big(|\Omega| + \|I_0\|^2_{L^2}\big)\,ds\Big\}\,.
\end{align*}
Since $g(u)$ has a positive lower bound depending only on the $u_0, T, \xi, \Omega$, from the above ineuality we see that there exists a constant $C_T$, depending only on 
$I_0, G_{\xi}, T$ and $h$ such that \eqref{uniform_final} holds true. This completes the proof. 

\noindent{\bf Proof of Theorem \ref{thm:weak}:}
Thanks to Theorem \ref{thm:wellposedness-local-solution} and Lemma \ref{lem:uniform-apriori-estimate}, there exists a unique local solution $(I,u,v)$ of the underlying problem \eqref{maina}-\eqref{maind} over a time interval $(0,\T)$ for small $\T>0$, and satisfies 
the uniform moment estimate \eqref{uniform_final}. 
More precisely, we extend the solution on a sequence of intervals $(0,t_n]$ such that $t_n\rightarrow \T$. Then considering
the initial problem starting from $\T$, one can extend the solution up to a given final time $T$, thanks to \eqref{uniform_final}. This completes the proof. 
\section{Numerical Approximation}
\label{sec:numerical}
To solve the proposed model \eqref{maina}-\eqref{maind}, we construct an explicit finite difference scheme, which is taken as the most straightforward option in the literature. Let $\tau$ be the time step size and $\tilde{h}$ is the spatial step size. Let $\V^n_{i,j} $ denotes the approximate value of $ \V(x_i,y_j,t_n)$, where
$x_i=i\tilde{h}, \hspace{0.2cm}	i=0,1,2...,N$ $y_j=j\tilde{h},  \hspace{0.2cm}   j=0,1,2...,M$ $t_n=n\tau,\hspace{0.2cm}  \left(n=0,1,2,...\right)$
where $n$ indicates the number of iterations and $M \times N$ is the size of the image. We choose symmetric boundary conditions for all the three variables in the system \eqref{maina}-\eqref{maind}, given as follows:
$$\V_{-1,j}^n=\V_{0,j}^n,\,\,\V_{N+1,j}^n=\V_{N,j}^n,\hspace{0.2cm}	j=0,1,2,...M, \hspace{1.0cm} \V_{i,-1}^n=\V_{i,0}^n,\,\,\,\V_{i,M+1}^n=\V_{i,M}^n,\hspace{0.2cm}	i=0,1,2,...N.$$	
We use the following finite difference approximations to replace the derivative terms in the system \eqref{maina}-\eqref{maind} :
\begin{align*}
& \dfrac{\partial \V_{i,j}^n }{\partial t} \approx  \displaystyle\frac{\V_{i,j}^{n+1}-\V_{i,j}^n}{\tau}\,, \quad 
\dfrac{\partial^2 \V_{i,j}^n }{\partial t^2} \approx \displaystyle\frac{\V_{i,j}^{n+1}-2\V_{i,j}^n+\V_{i,j}^{n-1}}{\tau ^2}\,,\\
& \nabla_x \V_{i,j}^n \approx \displaystyle\frac{\V_{i+1,j}^n-\V_{i-1,j}^n}{2\tilde{h}}\,, \quad 
\nabla_y \V_{i,j}^n \approx \displaystyle\frac{\V_{i,j+1}^n-\V_{i,j-1}^n}{2\tilde{h}}\,, \\
& \nabla^{+}_x \V_{i,j}^n \approx \displaystyle\frac{\V_{i+1,j}^n-\V_{i,j}^n}{\tilde{h}}\,, \quad 
\nabla^{+}_y \V_{i,j}^n \approx \displaystyle\frac{\V_{i,j+1}^n-\V_{i,j}^n}{\tilde{h}}\,, \\
& \nabla^{-}_x \V_{i,j}^n \approx \displaystyle\frac{\V_{i,j}^n-\V_{i-1,j}^n}{\tilde{h}}\,, \quad 
\nabla^{-}_y \V_{i,j}^n \approx \displaystyle\frac{\V_{i,j}^n-\V_{i,j-1}^n}{\tilde{h}}\,, \\
& \Delta_x \V_{i,j}^n \approx  \frac{{\V_{i+1,j}^n-2\V_{i,j}^n+\V_{i-1,j}^n}}{\tilde{h}^2}\,, \quad 
\Delta_y \V_{i,j}^n \approx \frac{{\V_{i,j+1}^n-2\V_{i,j}^n+\V_{i,j-1}^n}}{\tilde{h}^2}\,, \\
& |\nabla \V_{i,j}^n| \approx \sqrt{(\nabla_x \V_{i,j}^n)^2 + (\nabla_y \V_{i,j}^n)^2}\,.
\end{align*}
The discrete form of the equation \eqref{maina} could be written as follows:
\begin{align}\label{discrete_system}
&(1+\alpha \tau)I_{i,j}^{n+1}=(2+\alpha \tau)I_{i,j}^n -I_{i,j}^{n-1}
+{\tau ^2} \left[ \nabla^{+}_x \left( g(u_{i,j}^{n+1}) \nabla^{-}_x I_{i,j}^n \right) 
 + \nabla^{+}_y \left( g(u_{i,j}^{n+1}) \nabla^{-}_y I_{i,j}^n \right) \right] \nonumber \\
& \hspace{10.0cm}-2\tau ^2{\lambda^n} v^{n+1}_{i,j},\,\,\,\, n=1,2,3,...
\end{align}
with the initial condition $ I_{i,j}^0 =I_0(x_i,y_j)\,, \quad  I_{i,j}^1=I_{i,j}^0,\,\,\,\,	0 \leq i \leq N\,, ~~ 0 \leq j \leq M\,,$ $ \lambda^n $ is calculated using the formula as mentioned in $\left[ (2.9\text{c}) \cite{rudin1992nonlinear} \right]$.  $u_{i,j}^{n+1}$ and $v_{i,j}^{n+1}$ are calculated from the discretized equations of \eqref{mainb} and \eqref{mainc} respectively, given as follows
\begin{align}\label{discrete_u}
&(1+\beta \tau)u_{i,j}^{n+1}=(2+\beta\tau)u_{i,j}^n -u_{i,j}^{n-1}
+\kappa {\tau ^2} \left[ h_{i,j}^n-u_{i,j}^n+\frac{\nu^2}{2}\Delta u_{i,j}^n \right],\,\,\,\, n=1,2,3,...
\end{align}
and 
\begin{align}\label{disc:v}
&v_{i,j}^{n+1}=v_{i,j}^{n} + \tau \big(\Delta_x v_{i,j}^{n} +  \Delta_y v_{i,j}^{n}\big)-\tau\left(I_{i,j}^{0}-I_{i,j}^{n}  \right) \,,\,\,\, n=0,1,2,... 
\end{align}
with the initial conditions ${u}_{i,j}^0=G_\xi \ast |\nabla I_0|^{2},\,\,\,  u_{i,j}^1=u_{i,j}^0$, ${v}_{i,j}^0=0,$ and $ h_{i,j}^n= h_{\xi}\big(|\nabla(G_{\xi} \ast I^n_{i,j})| \big)$. We choose the function $h$ as $h(\theta)= 0.1 +\text{min}\{ \theta^2, K \}$ for numerical experiments, where $K$ is square of the maximum gray level value of the image $I$. Apart from the numerical discretization of \eqref{maina}-\eqref{maind}, a convergence criterion is required to stop the iterative process. To reach our destination, we start with a corrupted image $I_0$ and used the system \eqref{discrete_system}-\eqref{disc:v} repeatedly, resulting in a family of restored images ${I^p}$, which drafts the restored form of $I_0$. After satisfying the following condition
\begin{equation}\label{eq:4stopping}
\frac{{||I^{p+1}-I^p||^2_2}}{{||I^p||^2_2}}\leq \varepsilon,
\end{equation} 
we stop the iterative procedure. In \eqref{eq:4stopping} $I^p$ and $I^{(p+1)}$ illustrate  the restored images at the $p^{th}$ and ${(p+1)}^{th}$ iteration, respectively and $\varepsilon > 0$ is a predefined threshold. For the numerical experiments, we have used $\varepsilon = 10^{-4}$.

\section{Results and Discussion}
\label{sec:results}
\begin{figure}
       \begin{center}
        \begin{subfigure}[b]{0.3\textwidth}           
           \includegraphics[scale=0.25]{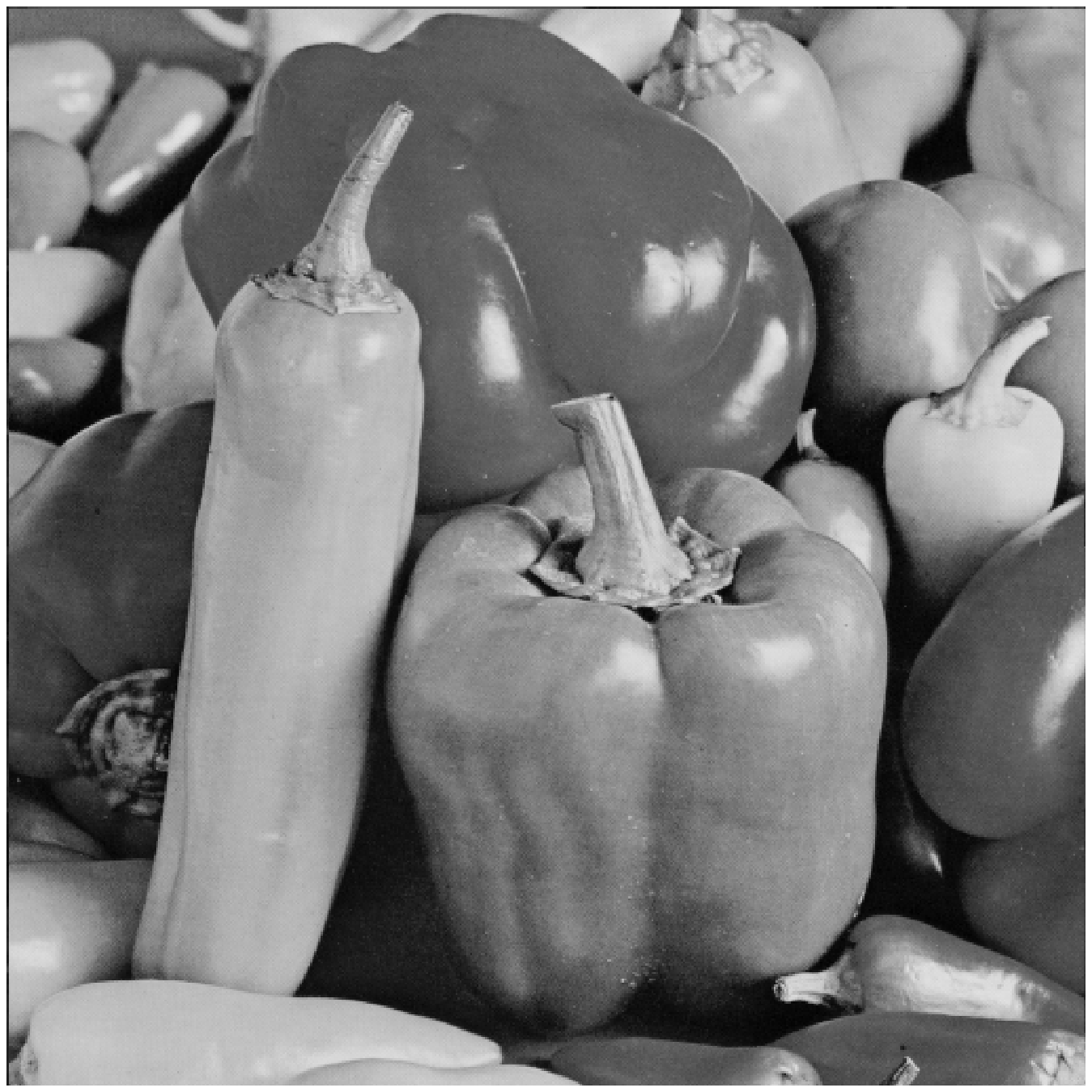}           
                \caption{Peppers}
                \label{fig:11a}
       \end{subfigure}
       \begin{subfigure}[b]{0.3\textwidth}           
                \includegraphics[scale=0.50]{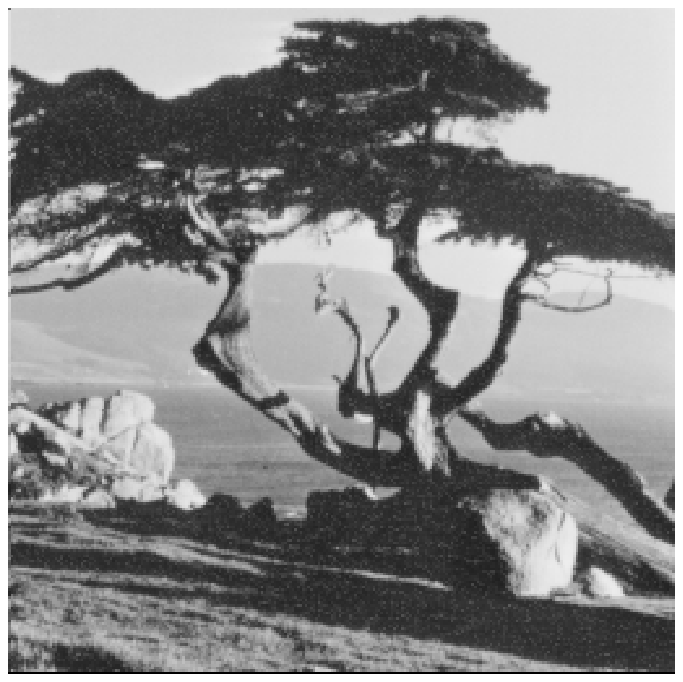}               
                \caption{Tree}
                \label{fig:11b}
       \end{subfigure}%
       \begin{subfigure}[b]{0.3\textwidth}           
                \includegraphics[scale=0.25]{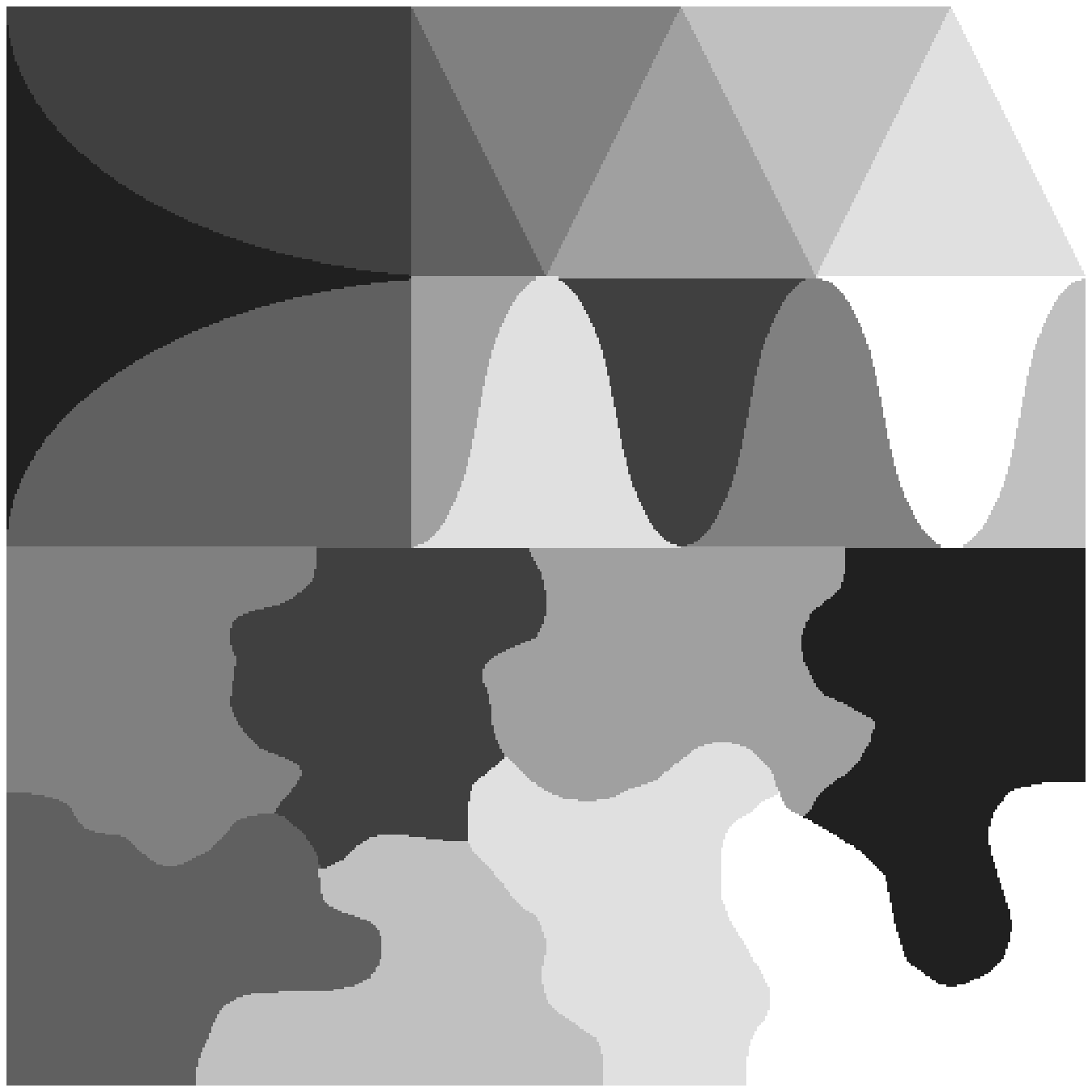}               
                \caption{Mosaic}
                \label{fig:11c}
       \end{subfigure}%
       
       \begin{subfigure}[b]{0.3\textwidth}           
                \includegraphics[scale=0.25]{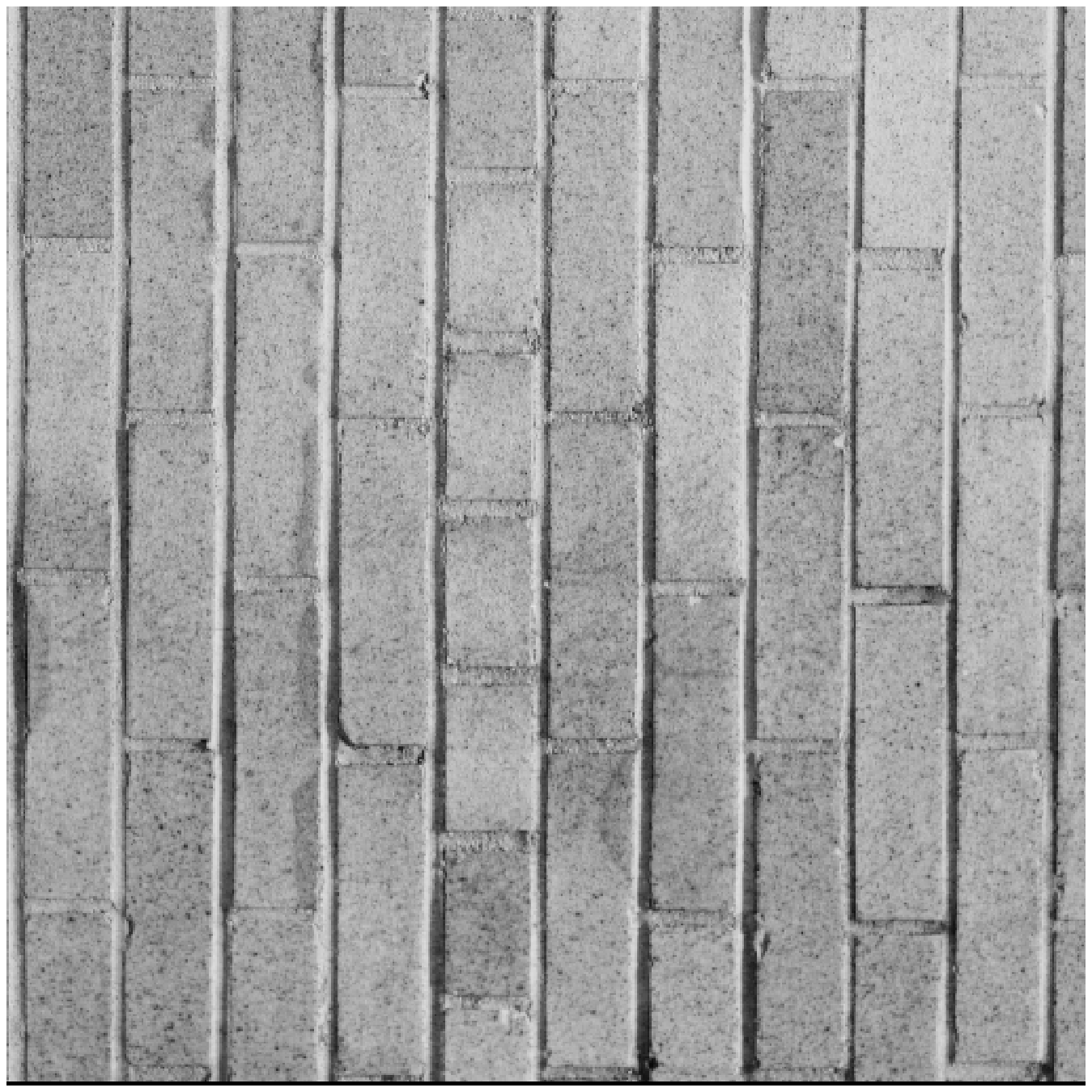}               
                \caption{Brick}
                \label{fig:11d}
       \end{subfigure}%
    \begin{subfigure}[b]{0.3\textwidth}           
                \includegraphics[scale=0.25]{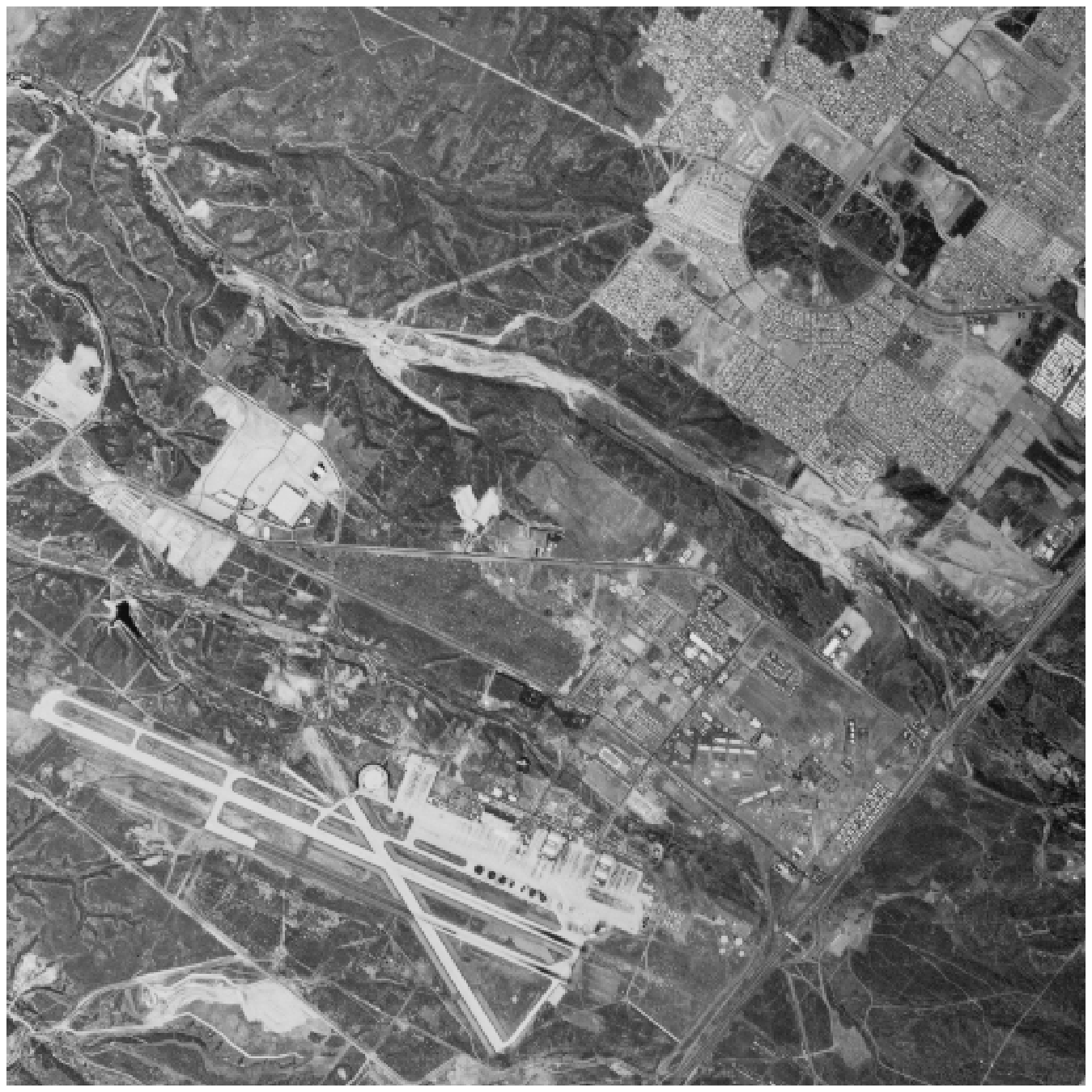}               
                \caption{Aerial}
                \label{fig:11e}
       \end{subfigure}%
      \begin{subfigure}[b]{0.3\textwidth}           
                \includegraphics[scale=0.50]{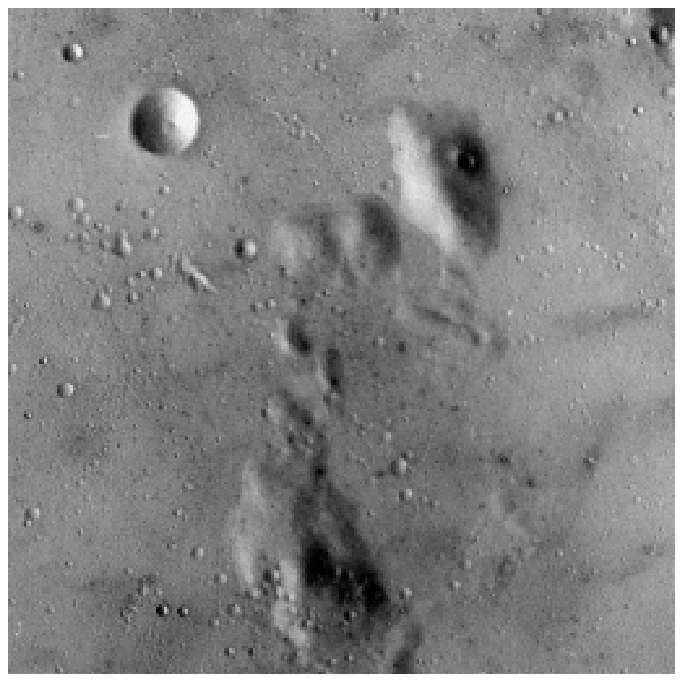}               
                \caption{Moon}
                \label{fig:11f}
       \end{subfigure}%
       
\caption{Set of test images: (a-b) natural images (c-d) texture images (e) aerial image (f) satellite image.}\label{fig:all images}
\end{center}
\end{figure}

Here, to assess the denoising performance of diffusion models, we present the experimental results of our proposed telegraph coupled partial differential equation model and compare against the results of TDE model \cite{ratner2007image}, Cao model \cite{cao2010class}, SYS model \cite{sun2016class}, and ACPDE model \cite{acpde2019rgate}. We tested all the filtering models on some standard gray level test images \ref{fig:all images}, which are degraded with additive Gaussian noise of zero mean($\mu$) and different level of standard deviations($\sigma$). The parameters used in the model are set to optimize the performance of the proposed algorithm. For all the methods, we choose time step size, $\tau=0.2$ and spatial step size, $\tilde{h}=1.$ Also, we choose $\xi=1$ to avoid over smoothing in the denoising process.

To objectively evaluate the denoising effect of all discussed algorithms, the performance is measured in terms of two commonly used quantitative metrics, Peak signal to noise ratio (PSNR)\cite{gonzalez2002digital} and mean structural similarity index measure (MSSIM)\cite{wang2009mean}. A higher values of MSSIM and PSNR illustrate the effective noise suppression.

\begin{figure}
       \centering

       \begin{subfigure}[b]{0.25\textwidth}           
                \includegraphics[scale=0.22]{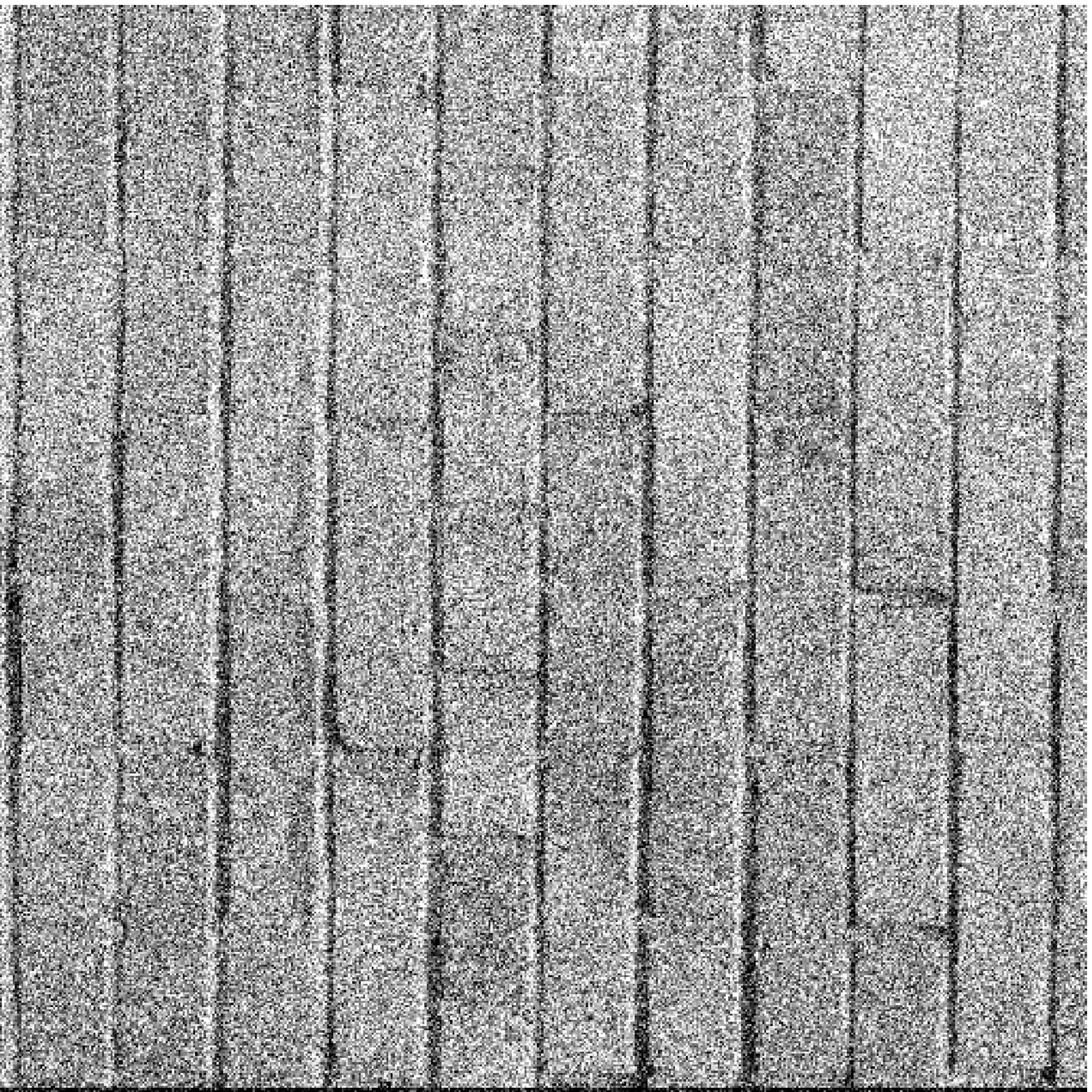}               
                \caption{Noisy}
                \label{fig:5b}
       \end{subfigure}%
       \begin{subfigure}[b]{0.25\textwidth}           
                \includegraphics[scale=0.22]{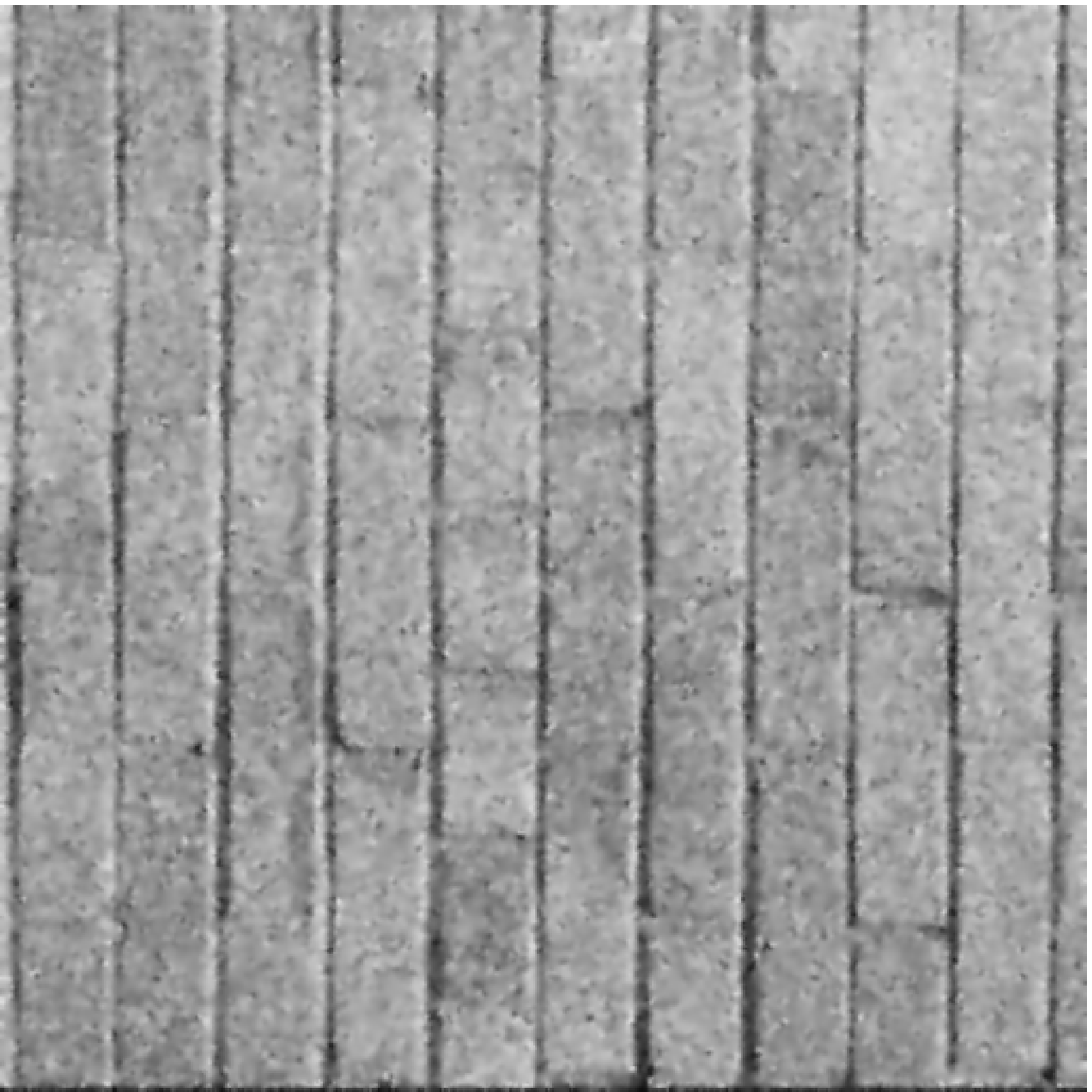}               
                \caption{SYS}
                \label{fig:5e}
       \end{subfigure}%
      \begin{subfigure}[b]{0.25\textwidth}           
                \includegraphics[scale=0.22]{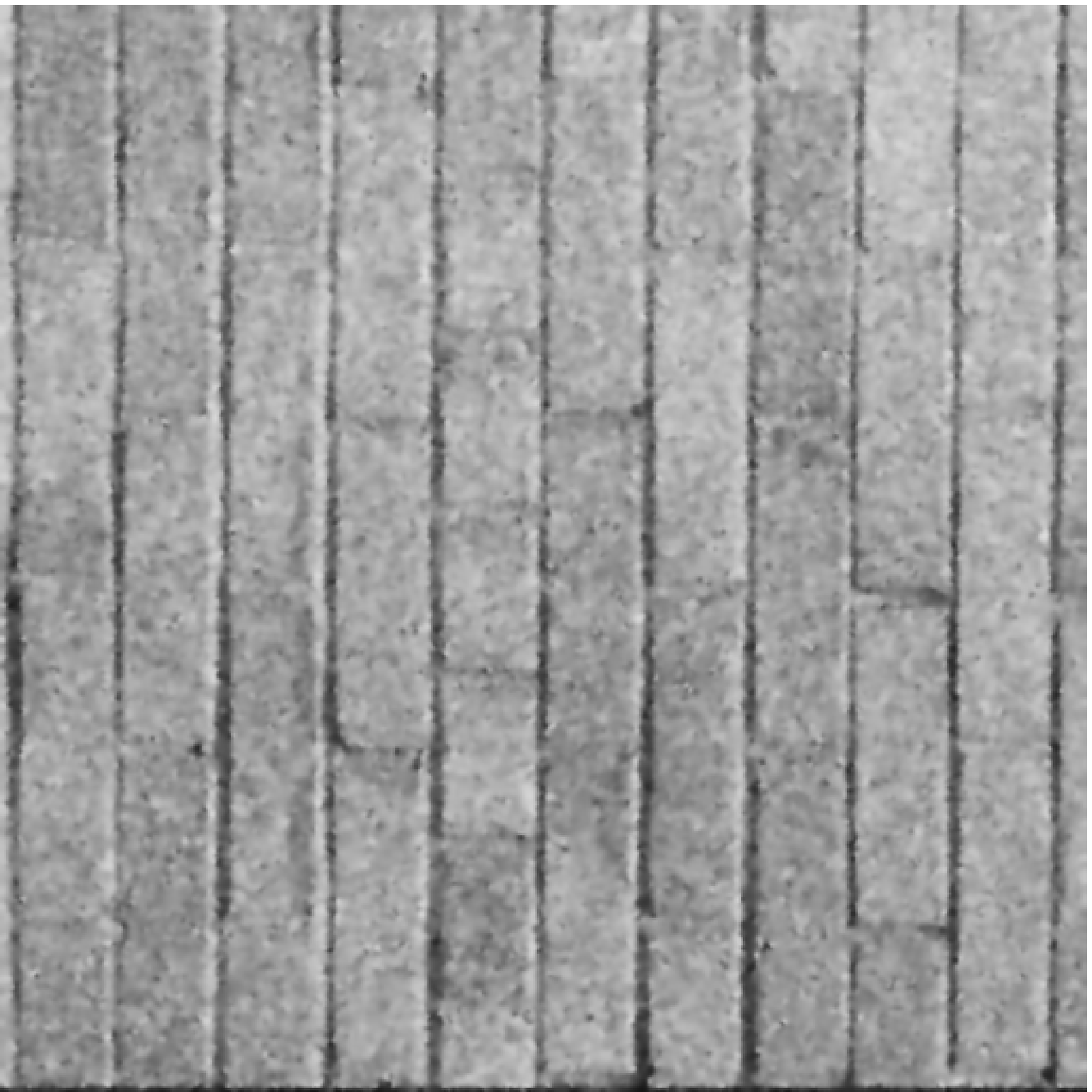}               
                \caption{ACPDE}
                \label{fig:5f}
       \end{subfigure}%
        \begin{subfigure}[b]{0.25\textwidth}           
                \includegraphics[scale=0.22]{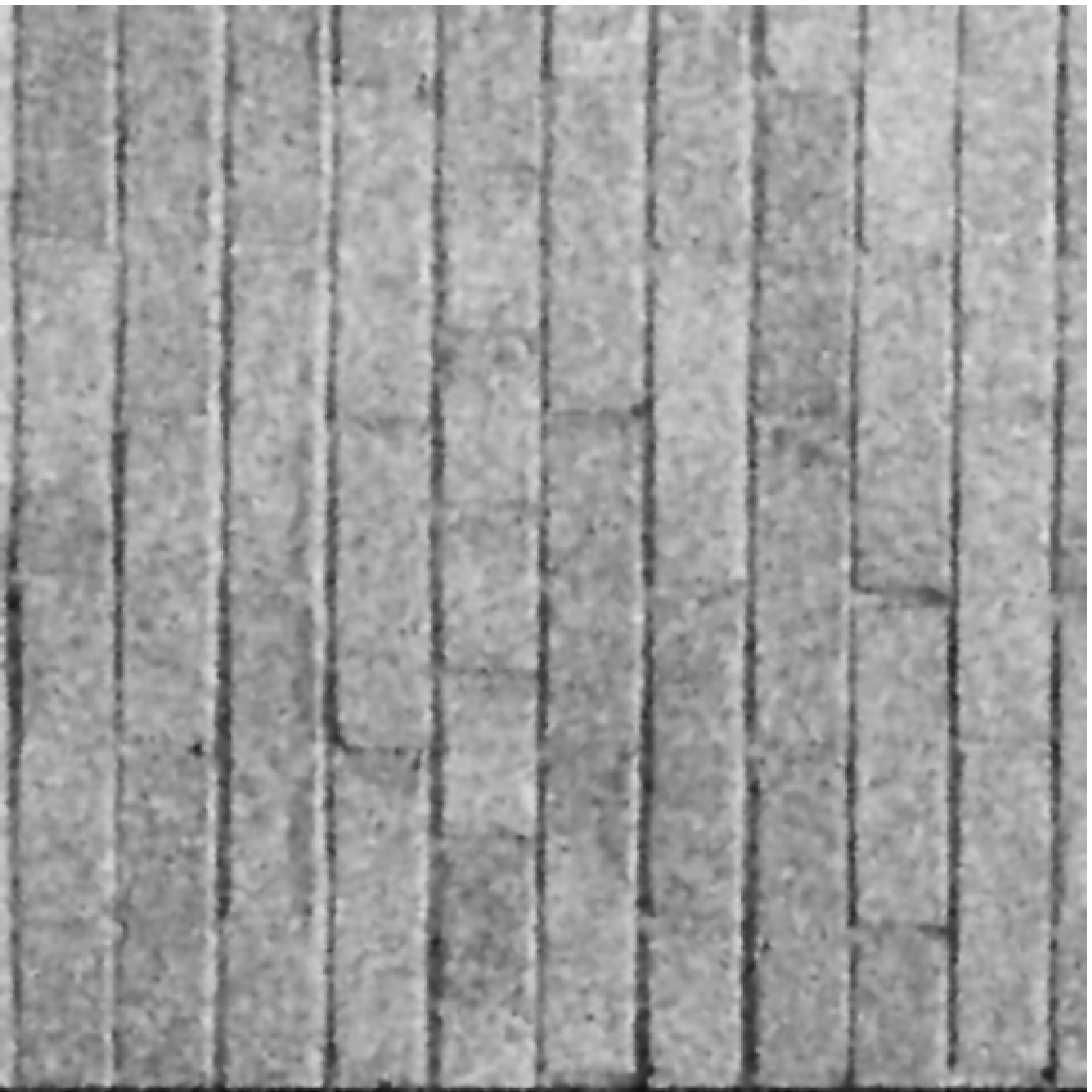}               
                \caption{Proposed}
                \label{fig:5g}
       \end{subfigure}%
       
 \caption{ A $512 \times 512$  brick image corrupted by Gaussian noise with $\sigma=60$ and restored by different models: (b) $\lambda =0.8, K=5$ (c) $\kappa=1, \nu=1, k=5$ (d) $\alpha=2, \beta=20,k=4.5, \nu=1.$}\label{fig:wall_60}
\end{figure}

Elaboration of results for a Brick image, corrupted with additive Gaussian noise of level $\sigma=60$, are shown in figure \ref{fig:wall_60}. This image contains a lot of fine texture. From the visual quality of denoised images, it is easy to perceive that the restored output obtained from the Proposed model gives better result in terms of denoising as well as edge-preserving. 

\begin{figure}
       \centering
       \begin{subfigure}[b]{0.3\textwidth}           
                \includegraphics[scale=0.28]{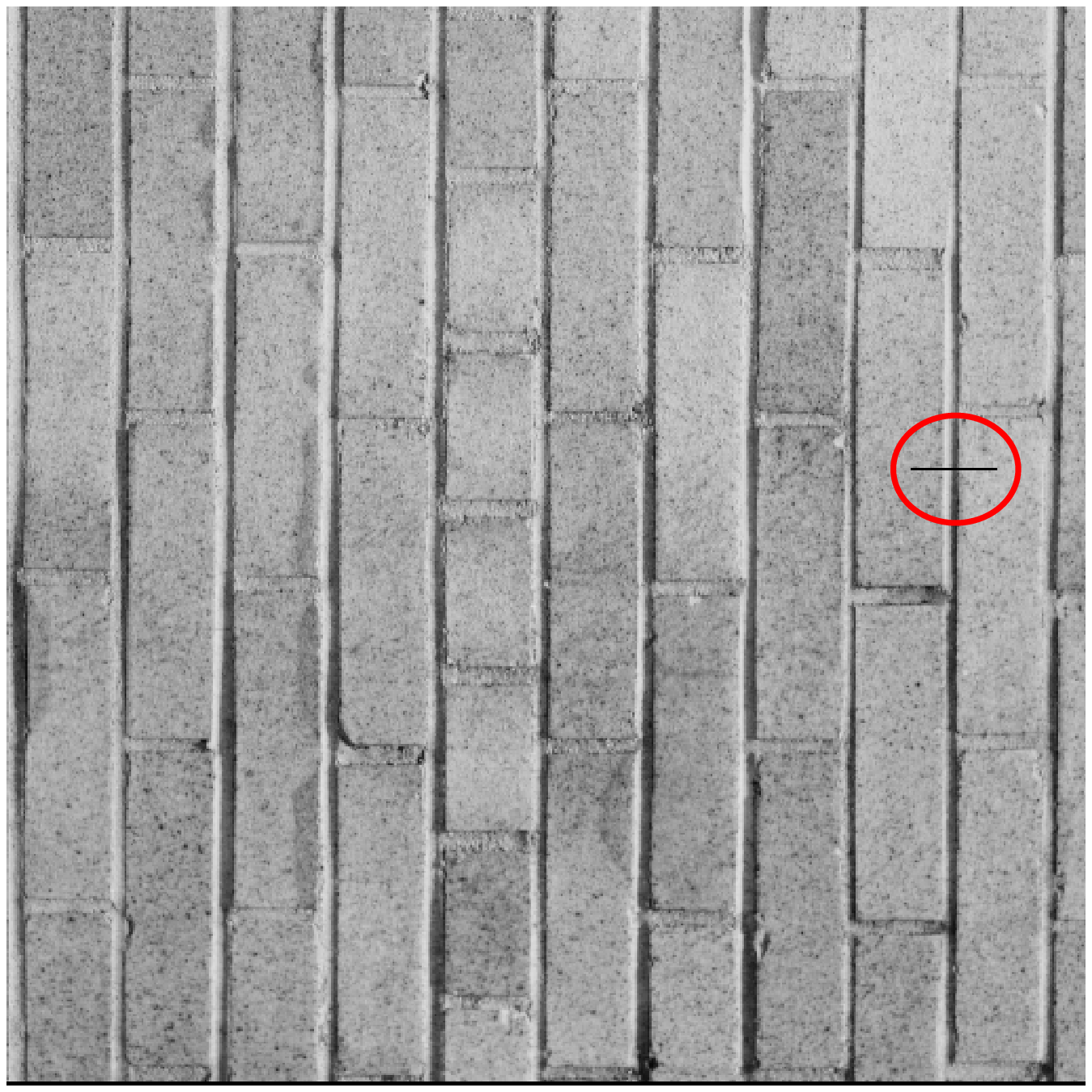}           
                \caption{}
                \label{fig:4a_wall_marked}
       \end{subfigure}%
              \begin{subfigure}[b]{0.4\textwidth}           
                \includegraphics[scale=0.4]{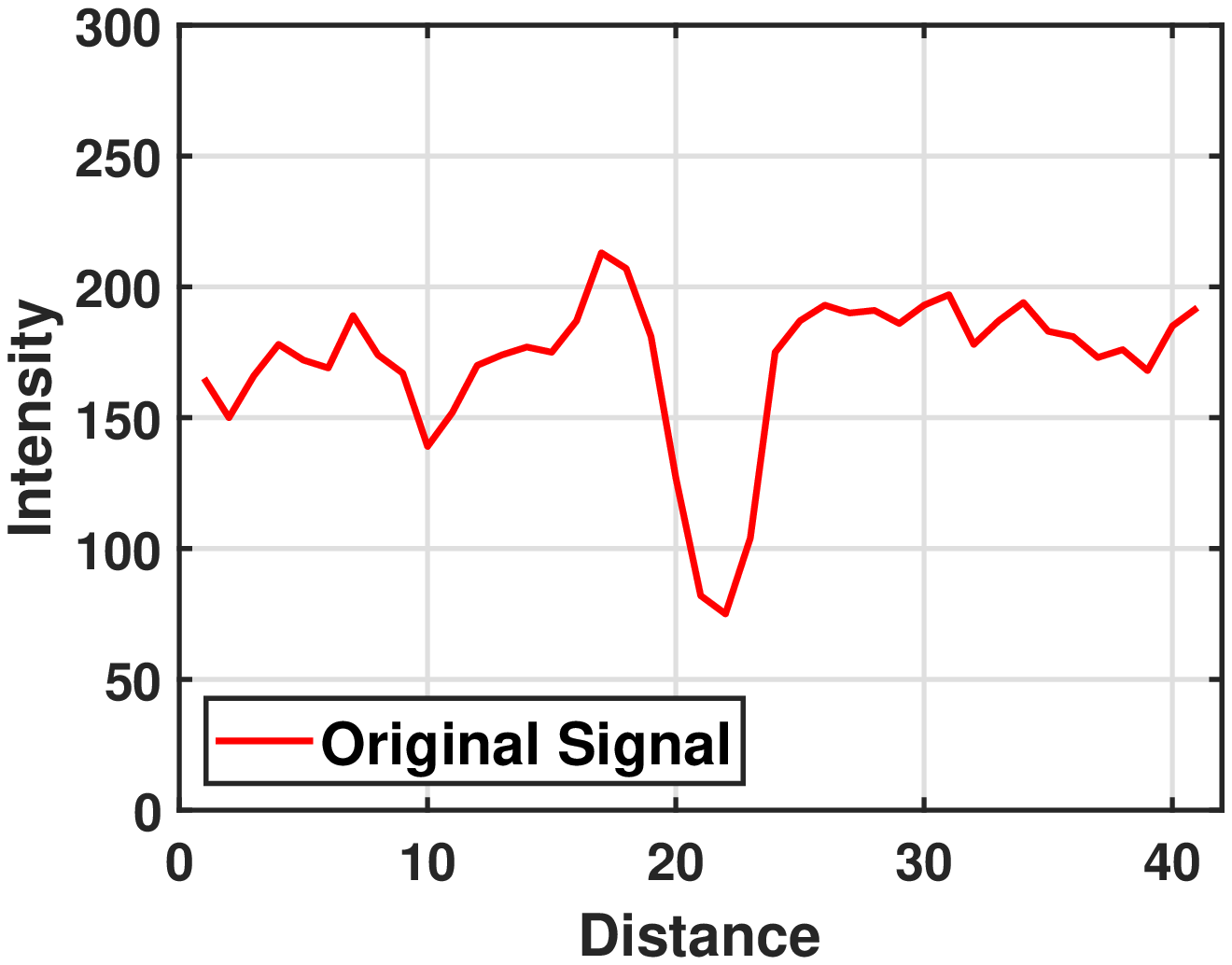}               
                \caption{}
                \label{fig:4a_wall_original}
       \end{subfigure}%

         \begin{subfigure}[b]{0.4\textwidth}           
                \includegraphics[scale=0.4]{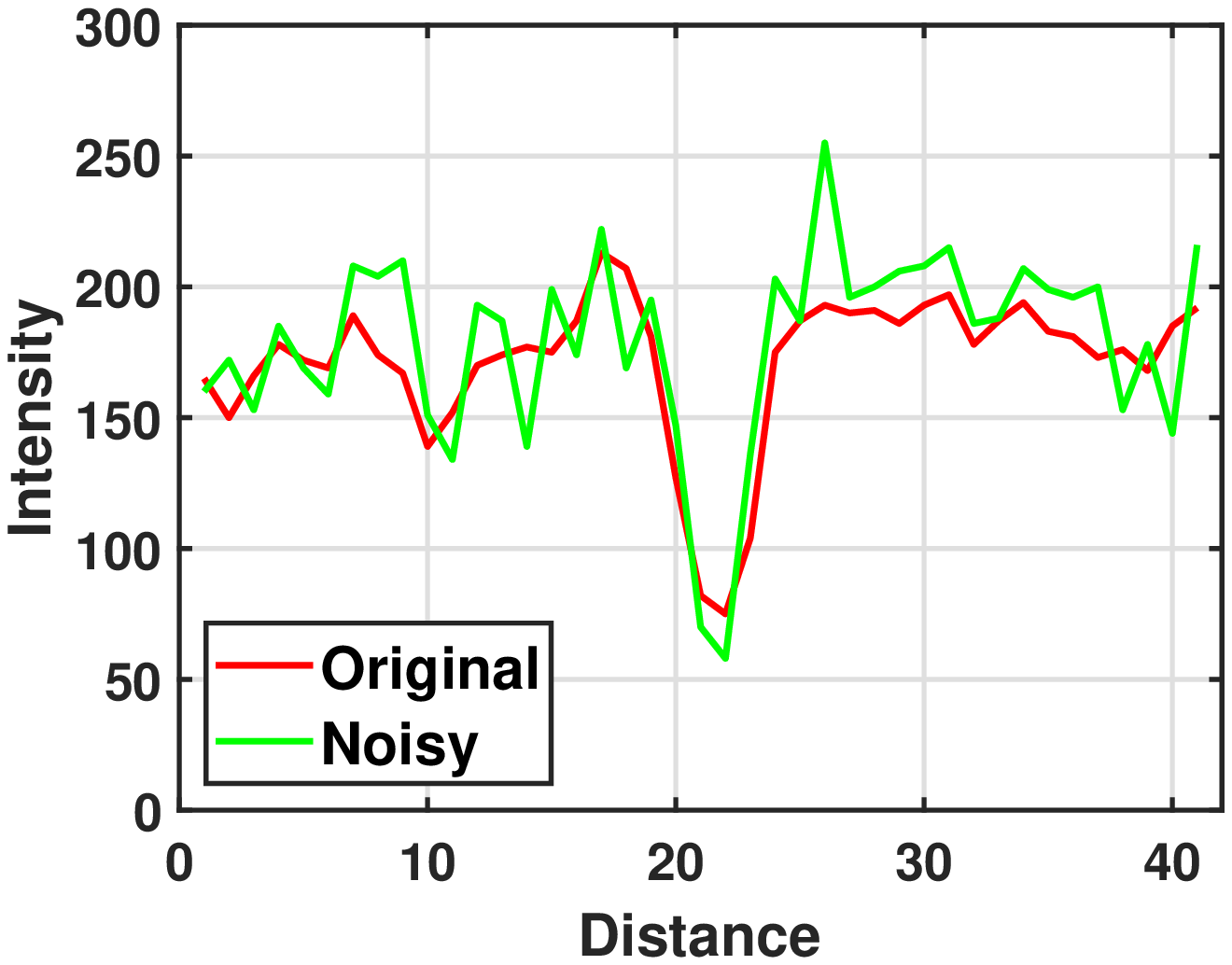}               
                \caption{}
                \label{fig:4a}
       \end{subfigure}%
         \begin{subfigure}[b]{0.4\textwidth}           
                \includegraphics[scale=0.4]{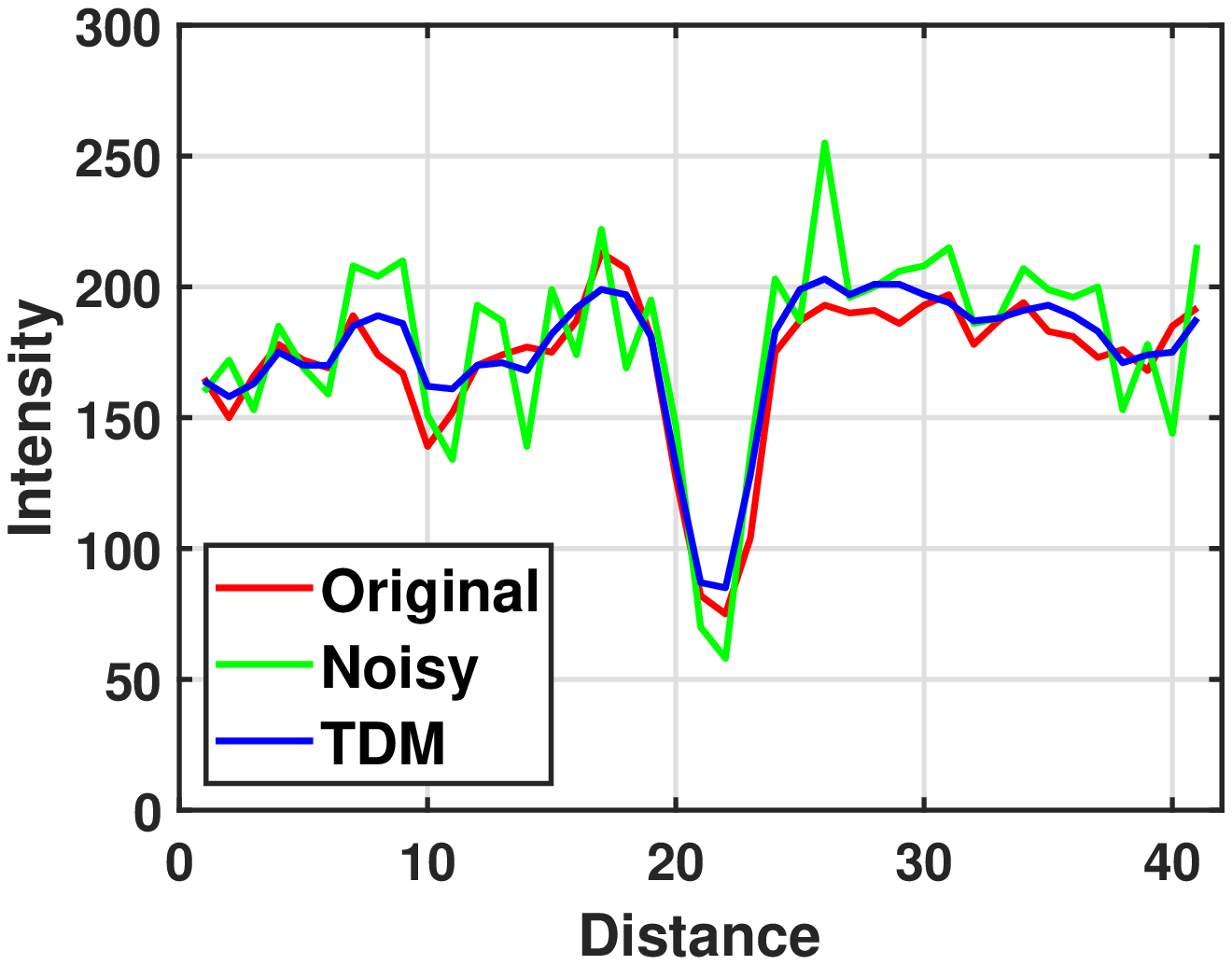}               
                \caption{}
                \label{fig:4b}
       \end{subfigure}%

       \begin{subfigure}[b]{0.4\textwidth}           
                \includegraphics[scale=0.4]{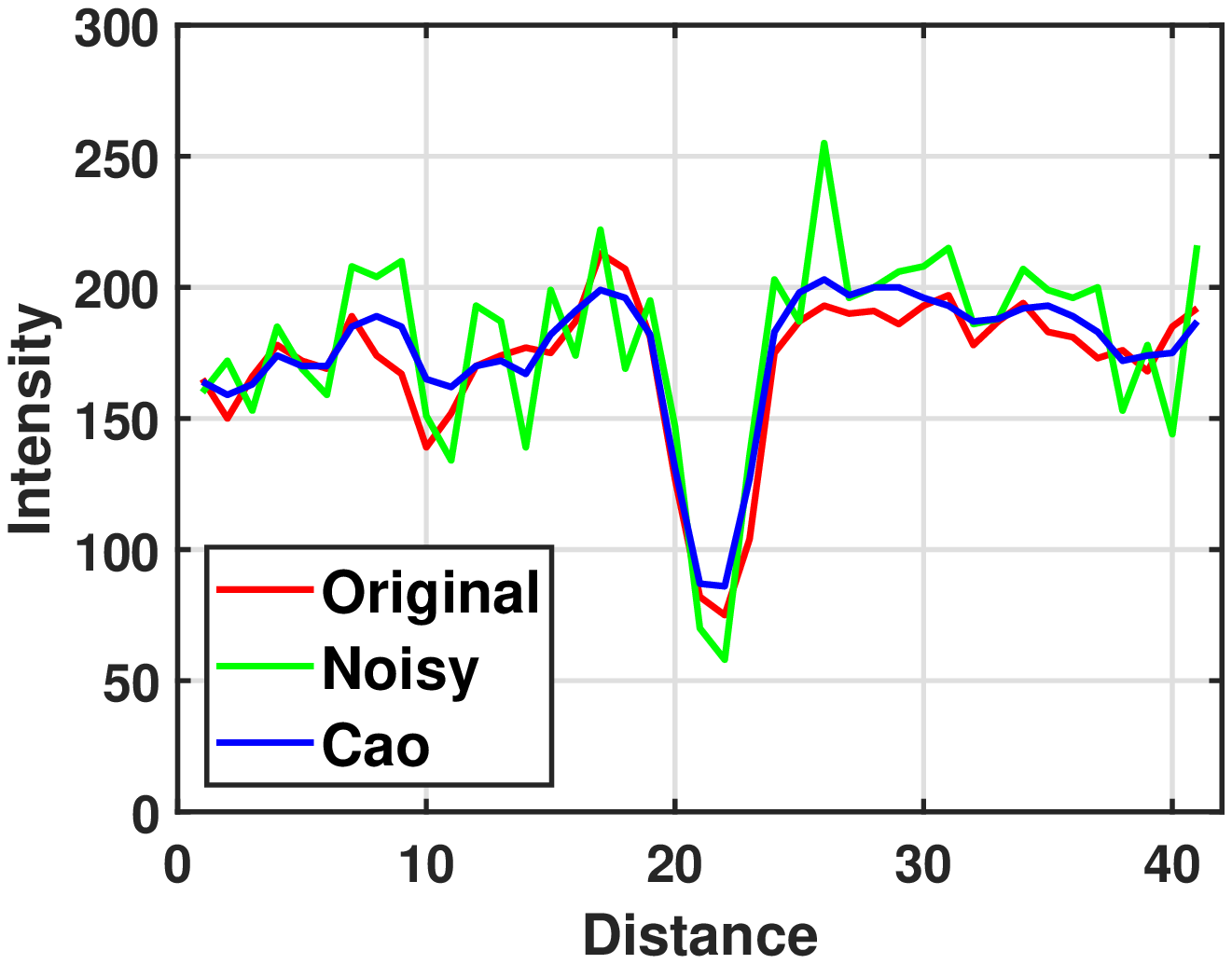}               
                \caption{}
                \label{fig:4c}
       \end{subfigure}%
      \begin{subfigure}[b]{0.4\textwidth}           
                \includegraphics[scale=0.4]{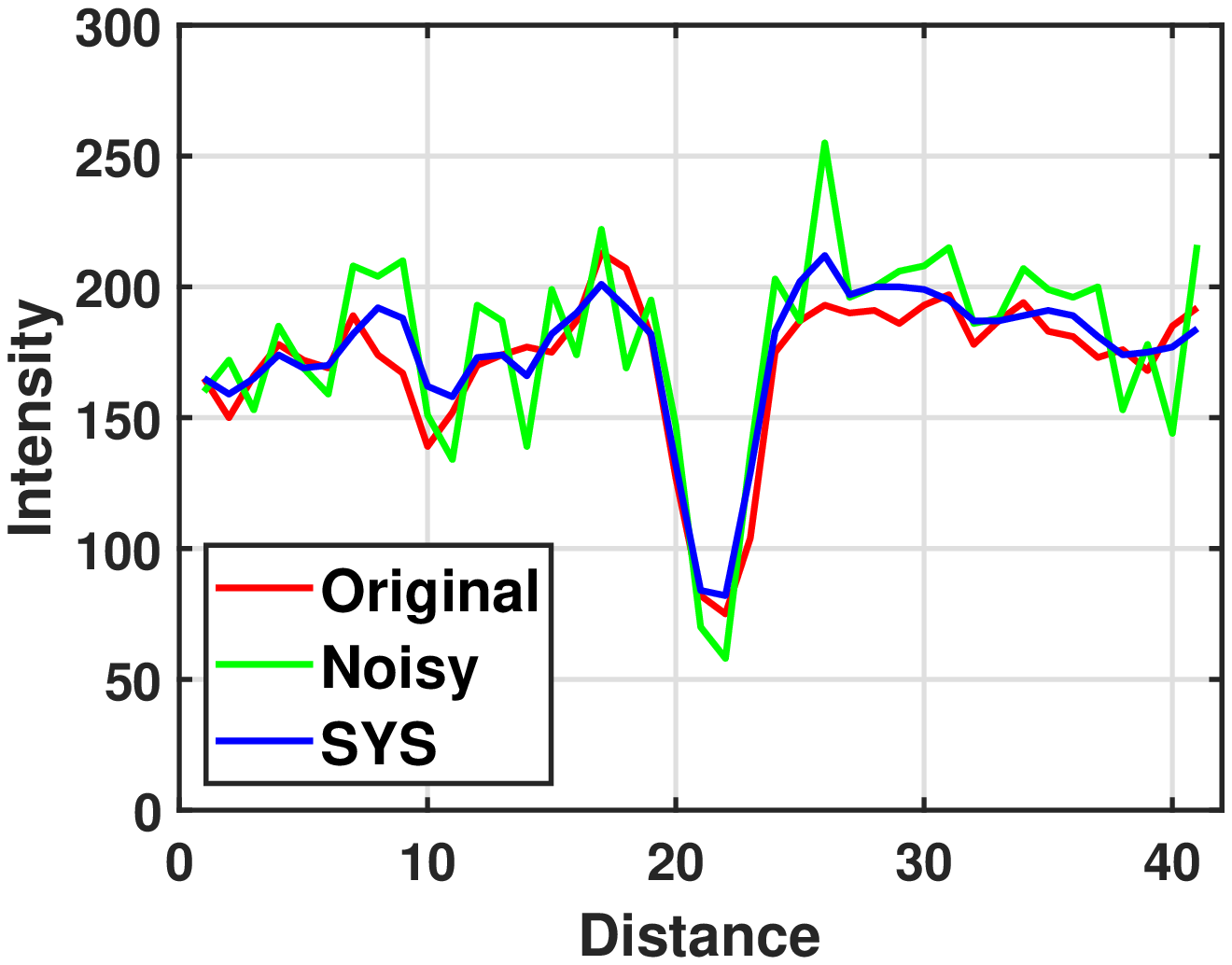}               
                \caption{}
                \label{fig:4d}
       \end{subfigure}%

         \begin{subfigure}[b]{0.4\textwidth}           
                \includegraphics[scale=0.4]{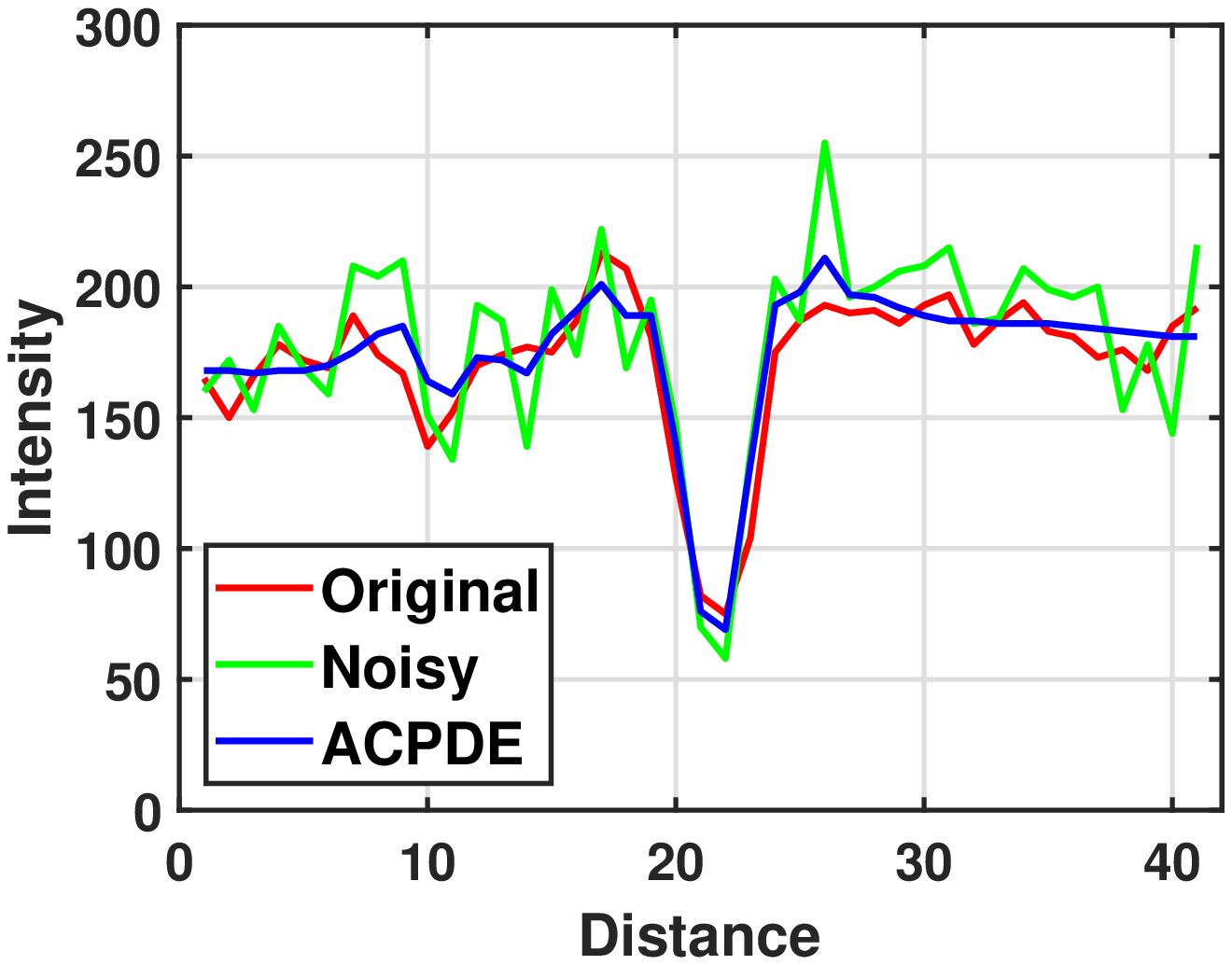}               
                \caption{}
                \label{fig:4e}
       \end{subfigure}%
        \begin{subfigure}[b]{0.4\textwidth}           
                \includegraphics[scale=0.4]{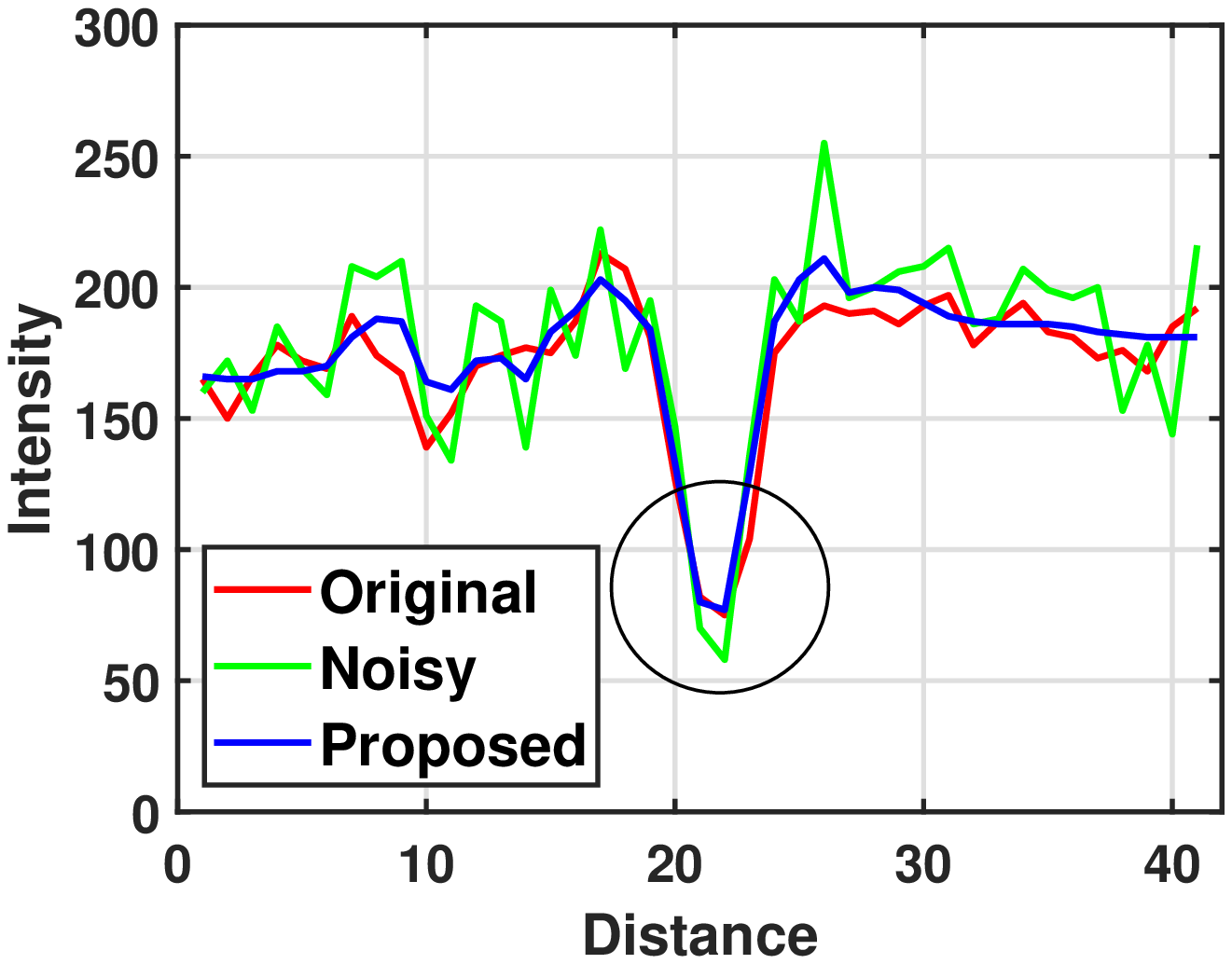}               
                \caption{}
                \label{fig:4f}
       \end{subfigure}%
       
 \caption{Restored signal by different models when the Brick image is corrupted by Gaussian noise with $\sigma=20$: (a) Original image showing the 1D slice (b) Clear signal (c) Noisy model (d) TDM; $ \lambda =5, K=15 $ (e)Cao; $\lambda =20, K=6$ (f) SYS; $\lambda =2, K=5$ (g) ACPDE; $\kappa=1, \nu=1, k=3$ (h) Proposed; $\alpha=1, \beta=20,k=1.15, \nu=1.$ }\label{fig:wall_20_signal}
\end{figure}

Along the qualitative analysis using full image surface, in figure \ref{fig:wall_20_signal}, we have also explored the quality of resultant images using a slice of a Brick image. This figure shows the slice of the original, noisy, and restored images. From figure \ref{fig:4f} it is easy to conclude that the restored signal obtained by the proposed model is more closure to the original signal in comparison with the other discussed models.
 
\begin{figure}
       \centering
     
       \begin{subfigure}[b]{0.22\textwidth}           
           \includegraphics[scale=0.19]{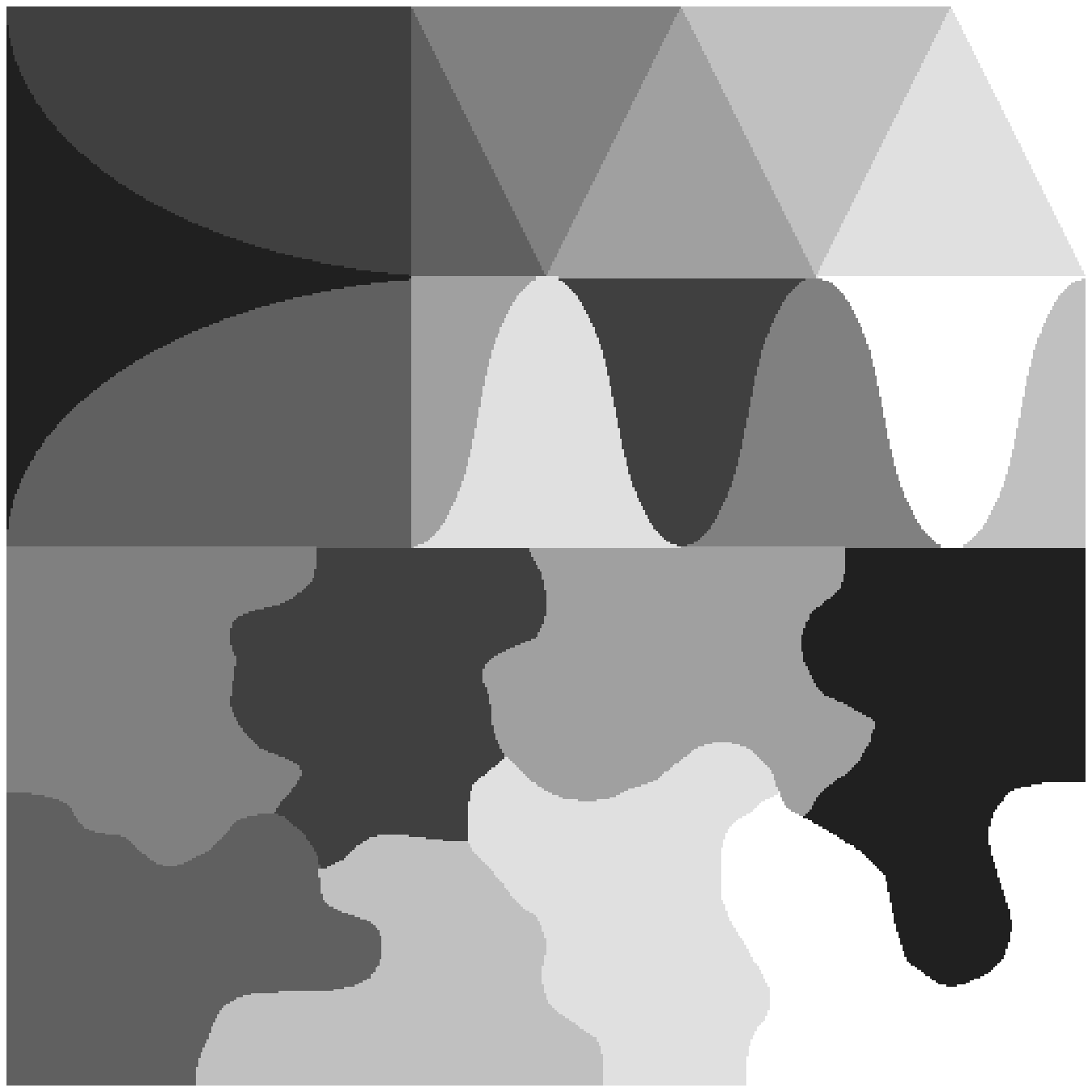}           
                \caption{}
                \label{fig:3a_clear}
       \end{subfigure}%
      \begin{subfigure}[b]{0.22\textwidth}           
                \includegraphics[scale=0.19]{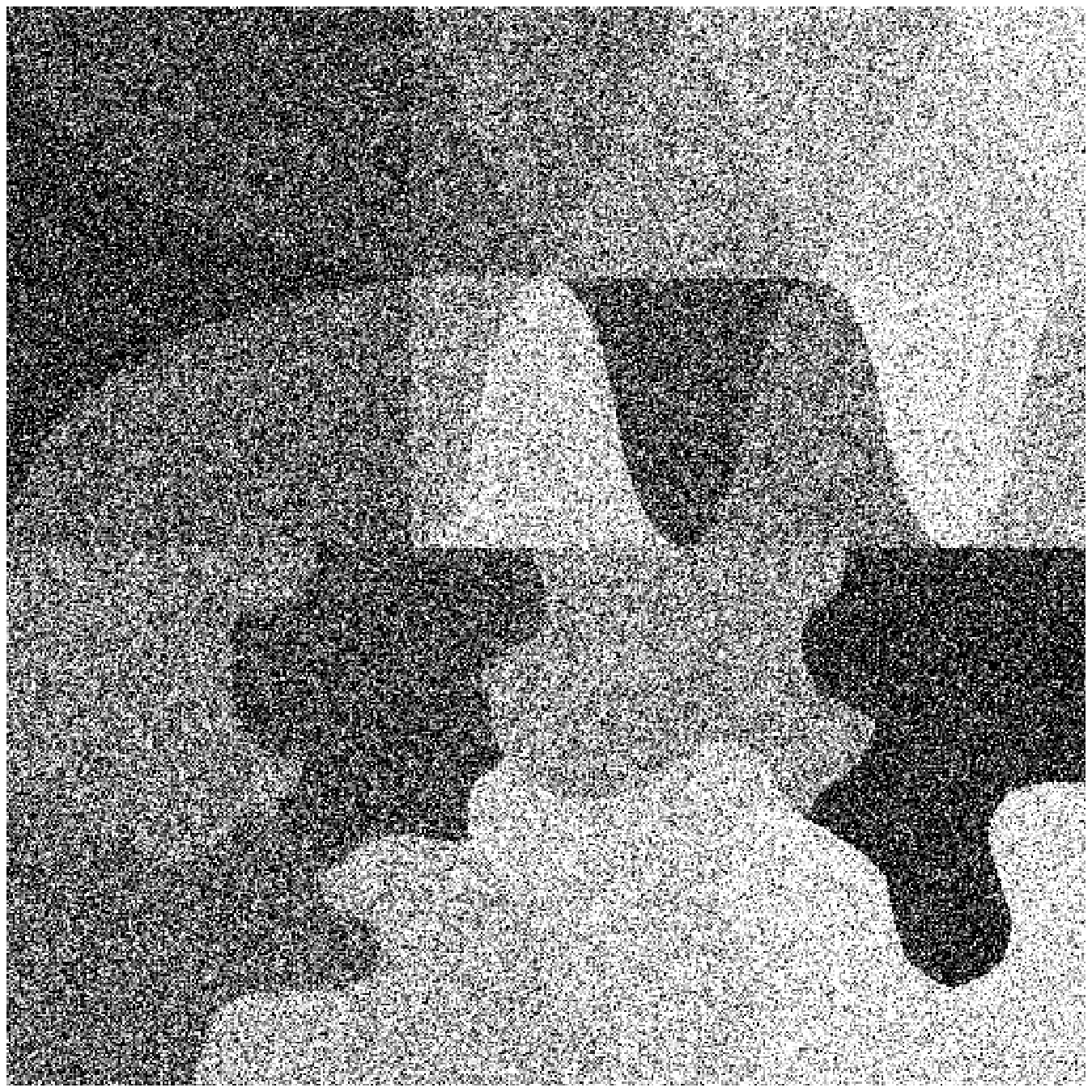}   
                    \caption{}         
                \label{fig:3a_noisy}
       \end{subfigure}%
           \begin{subfigure}[b]{0.22\textwidth}           
                \includegraphics[scale=0.22]{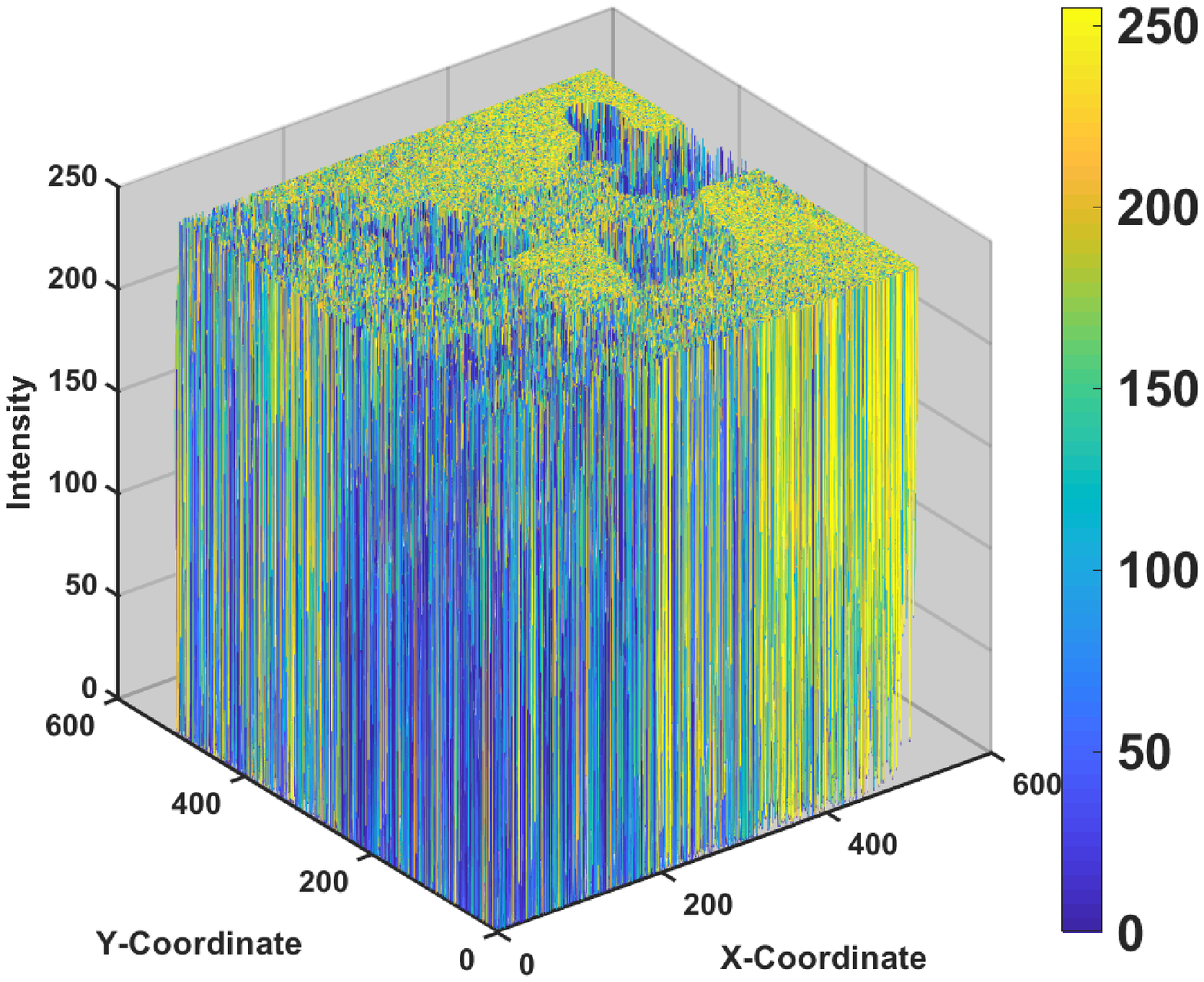}      
                     \caption{}     
                \label{fig:3a_noisy_surf}
       \end{subfigure}%

      \begin{subfigure}[b]{0.22\textwidth}           
                \includegraphics[scale=0.19]{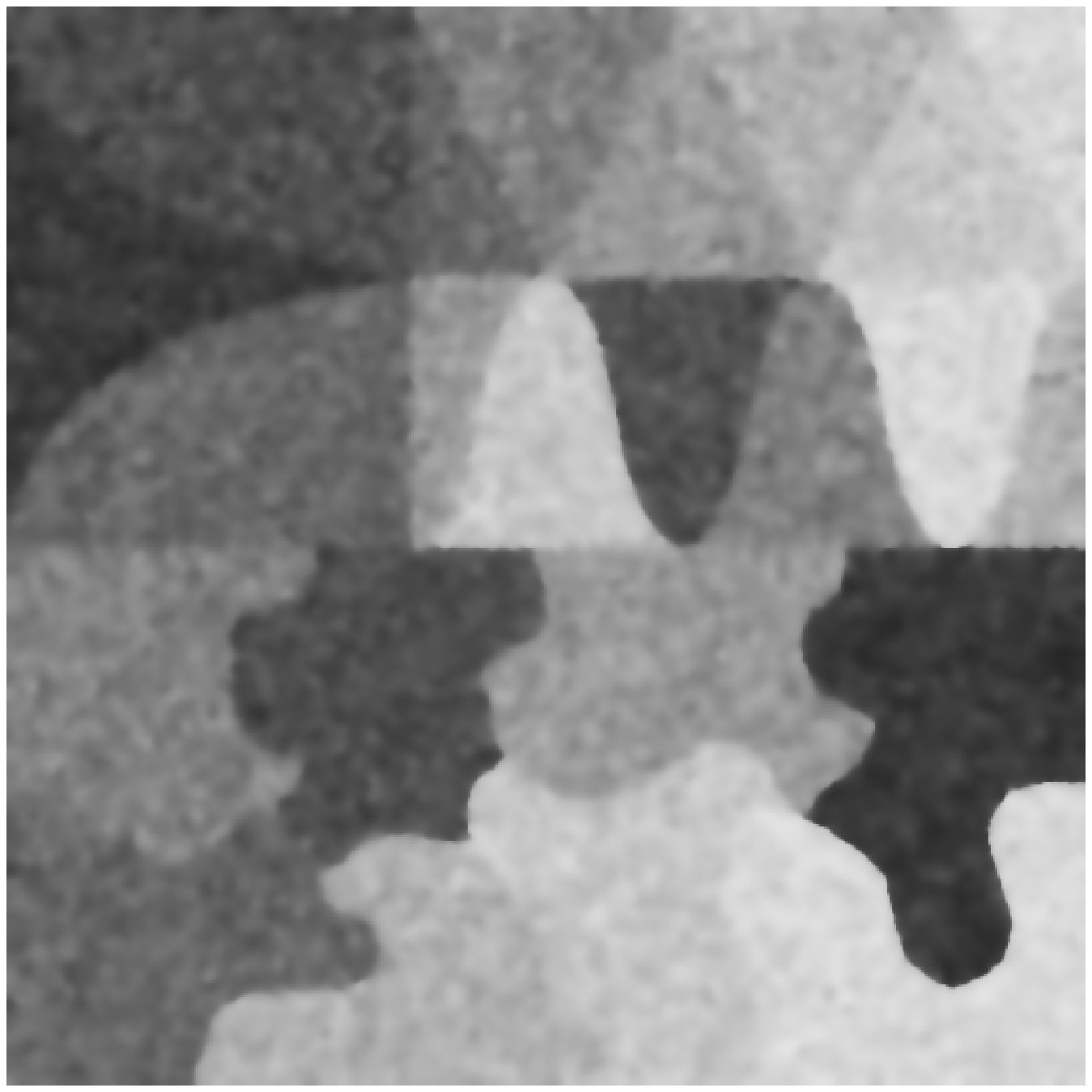}    
                     \caption{}       
                \label{fig:3c}
       \end{subfigure}%
              \begin{subfigure}[b]{0.22\textwidth}           
                \includegraphics[scale=0.19]{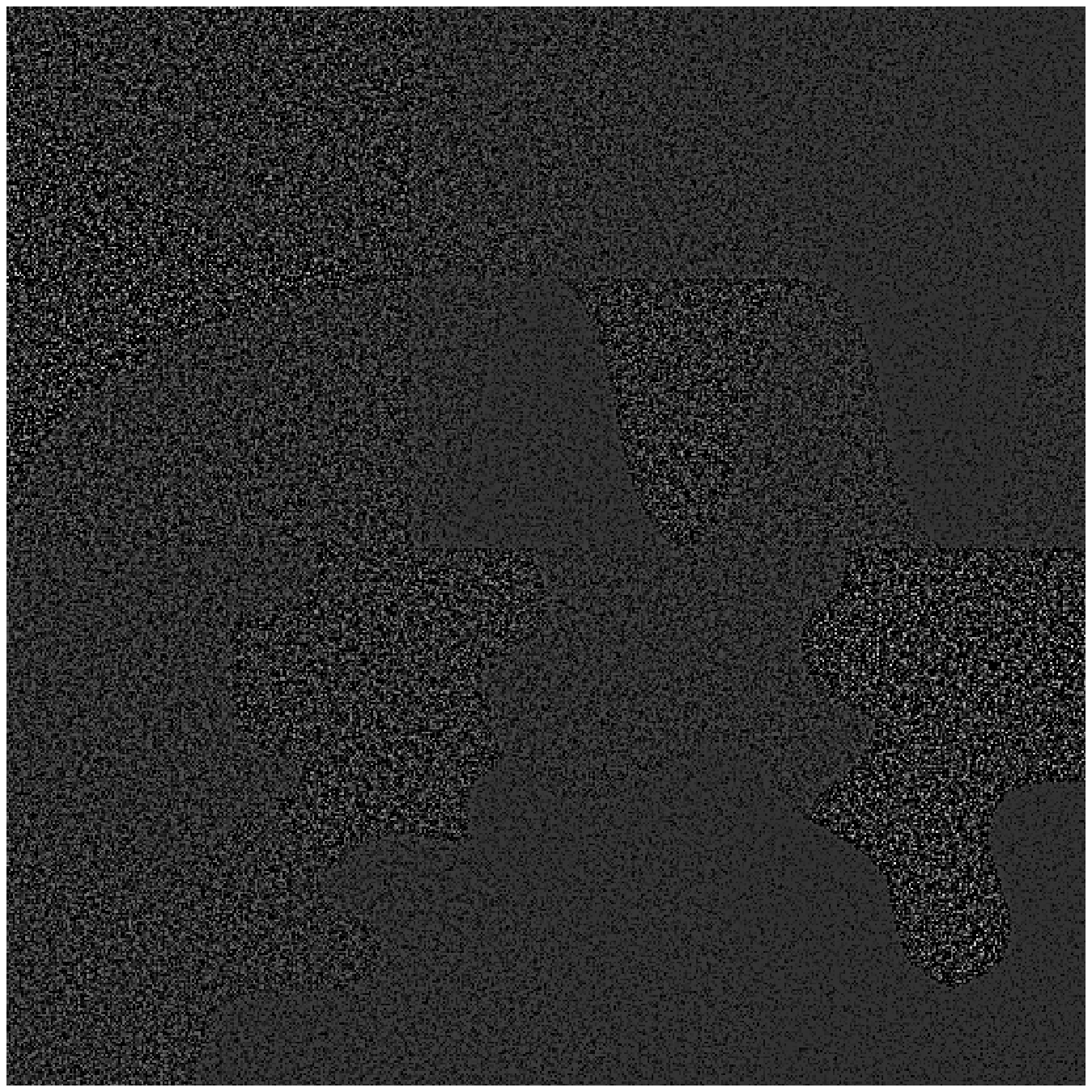} 
                     \caption{}          
                \label{fig:3c_ratio}
       \end{subfigure}%
                  \begin{subfigure}[b]{0.22\textwidth}           
                \includegraphics[scale=0.22]{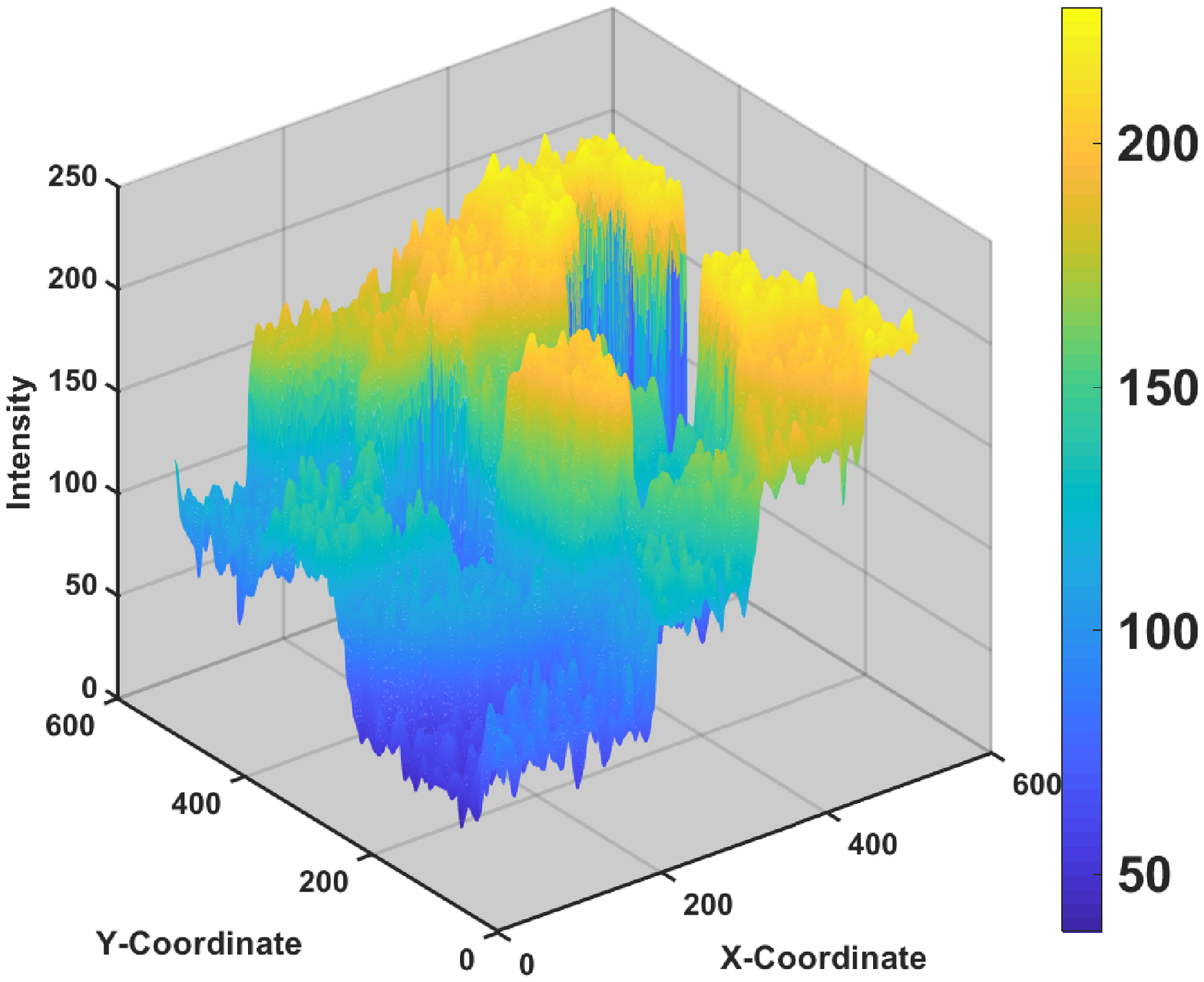} 
                       \caption{}        
                \label{fig:3c_surf}
       \end{subfigure}%

    \begin{subfigure}[b]{0.22\textwidth}           
                \includegraphics[scale=0.19]{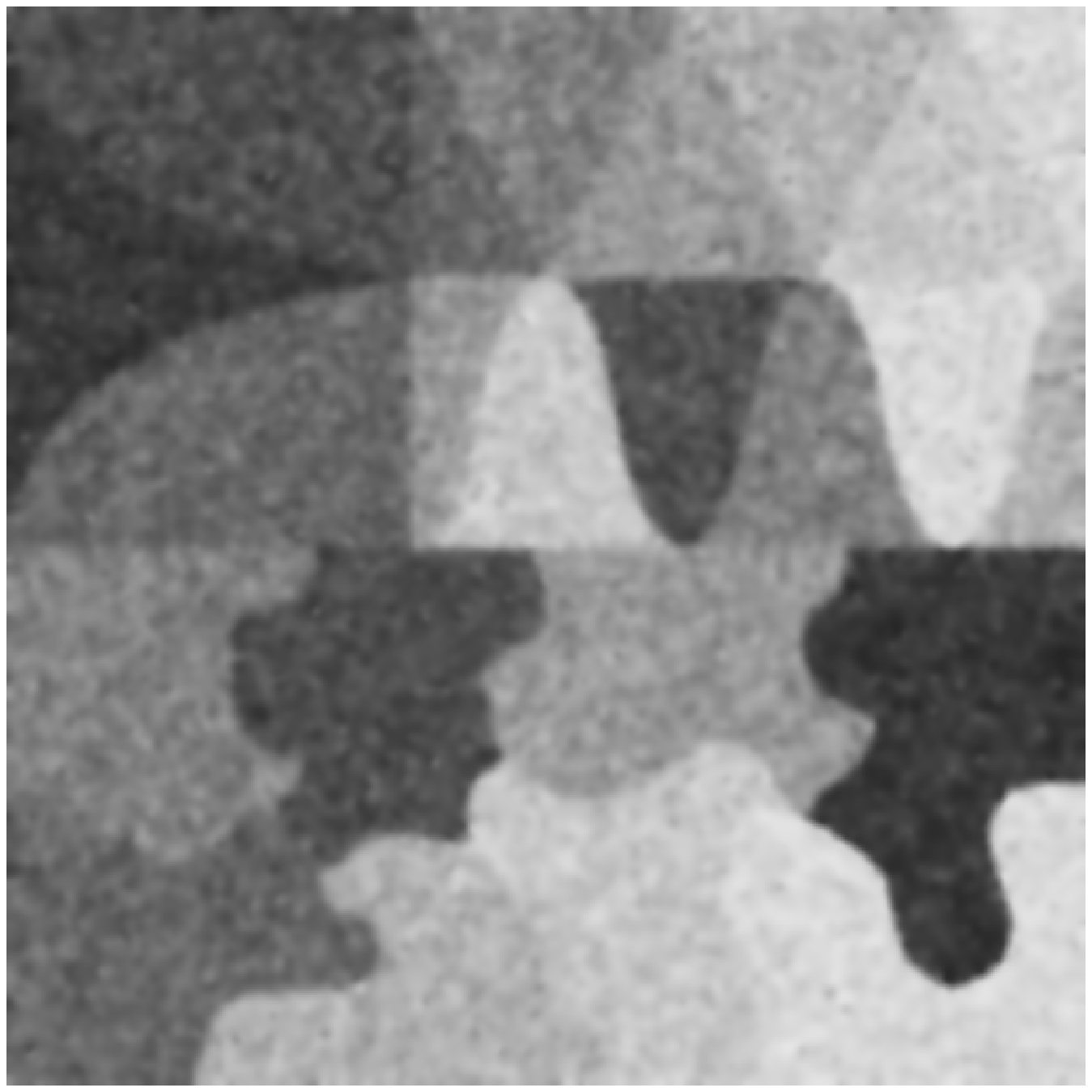}    
                       \caption{}     
                \label{fig:3d}
       \end{subfigure}%
           \begin{subfigure}[b]{0.22\textwidth}           
                \includegraphics[scale=0.19]{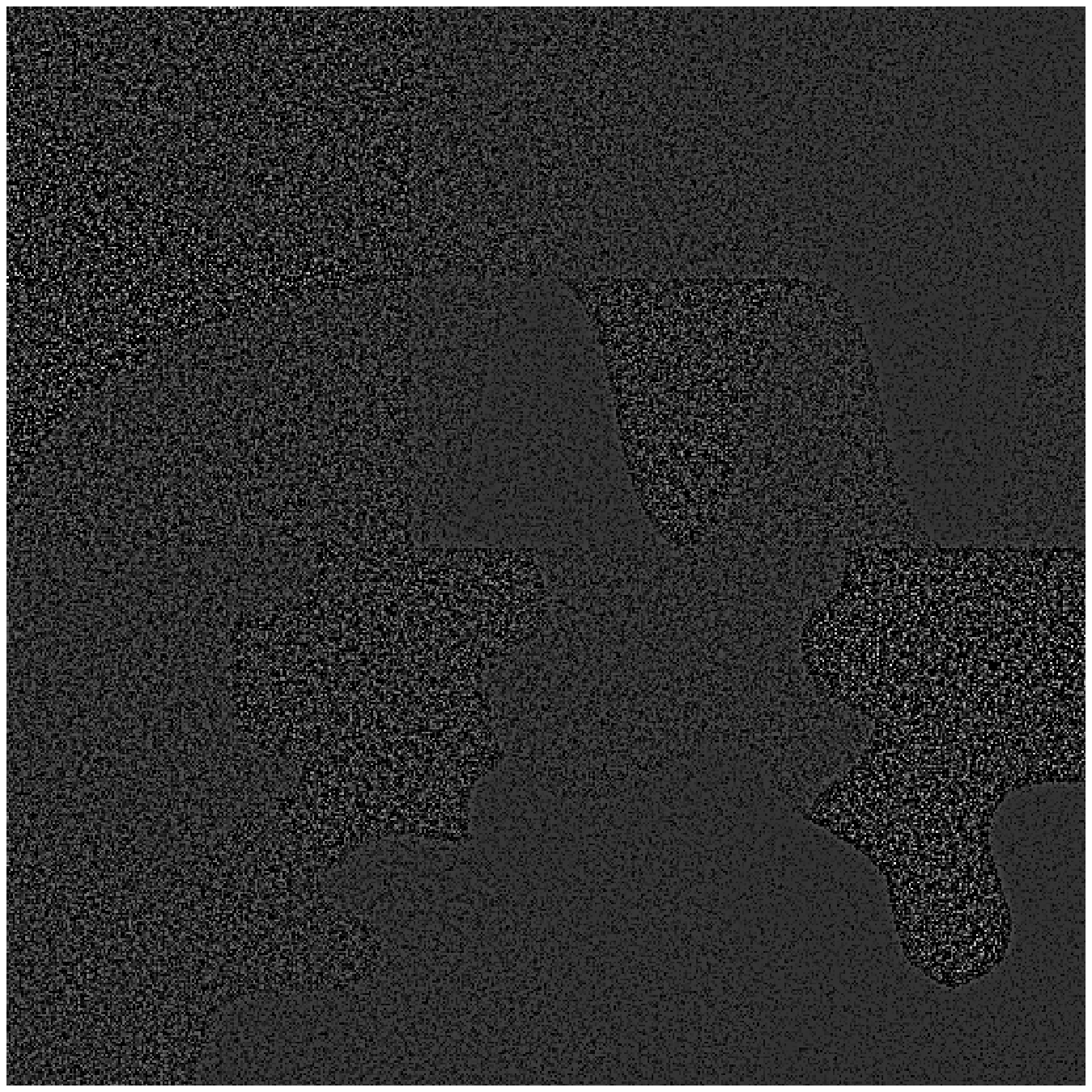}  
                        \caption{}      
                \label{fig:3d_ratio}
       \end{subfigure}%
                  \begin{subfigure}[b]{0.22\textwidth}           
                \includegraphics[scale=0.22]{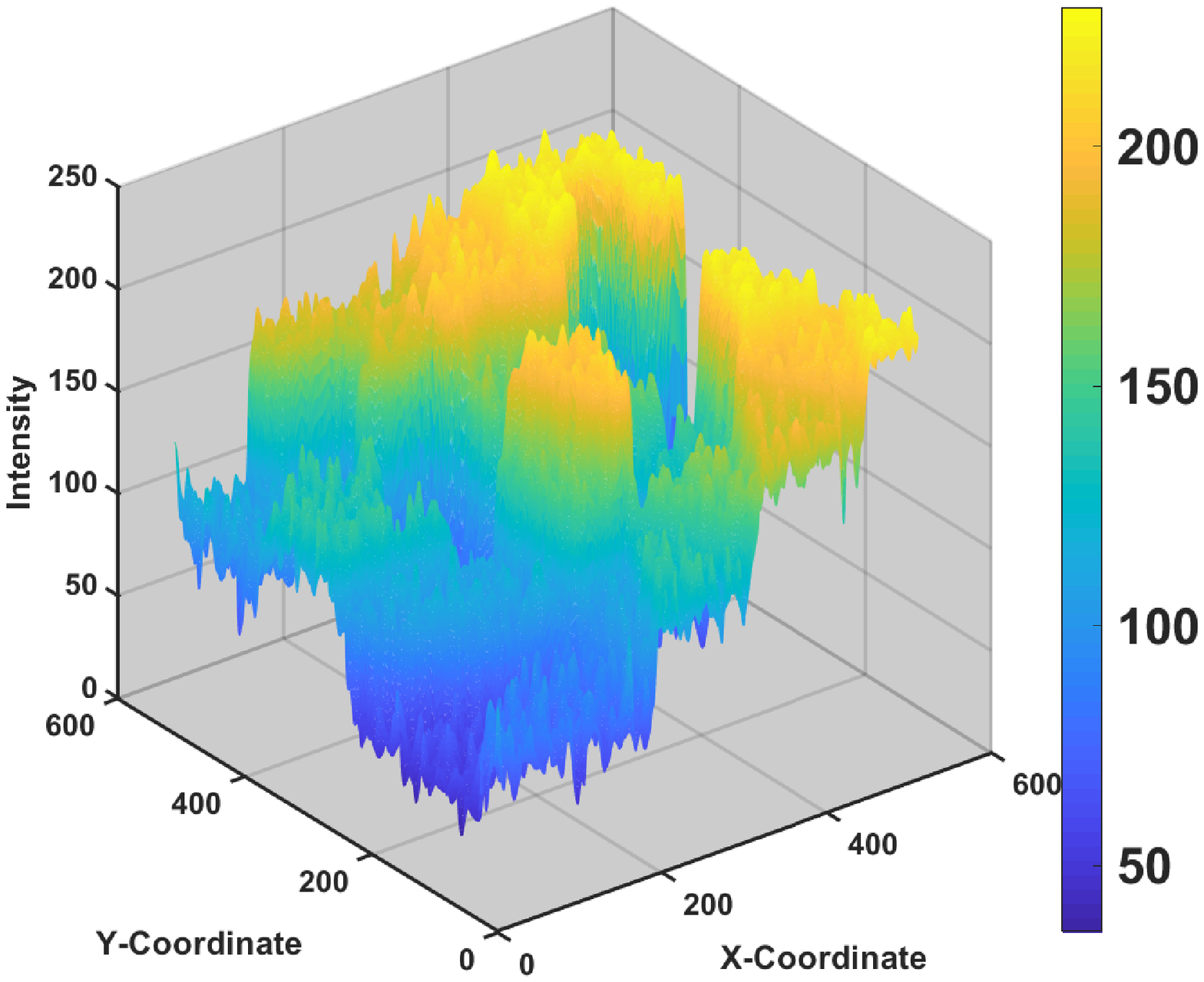}  
                         \caption{}     
                \label{fig:3d_surf}
       \end{subfigure}%
       
      \begin{subfigure}[b]{0.25\textwidth}           
                \includegraphics[scale=0.19]{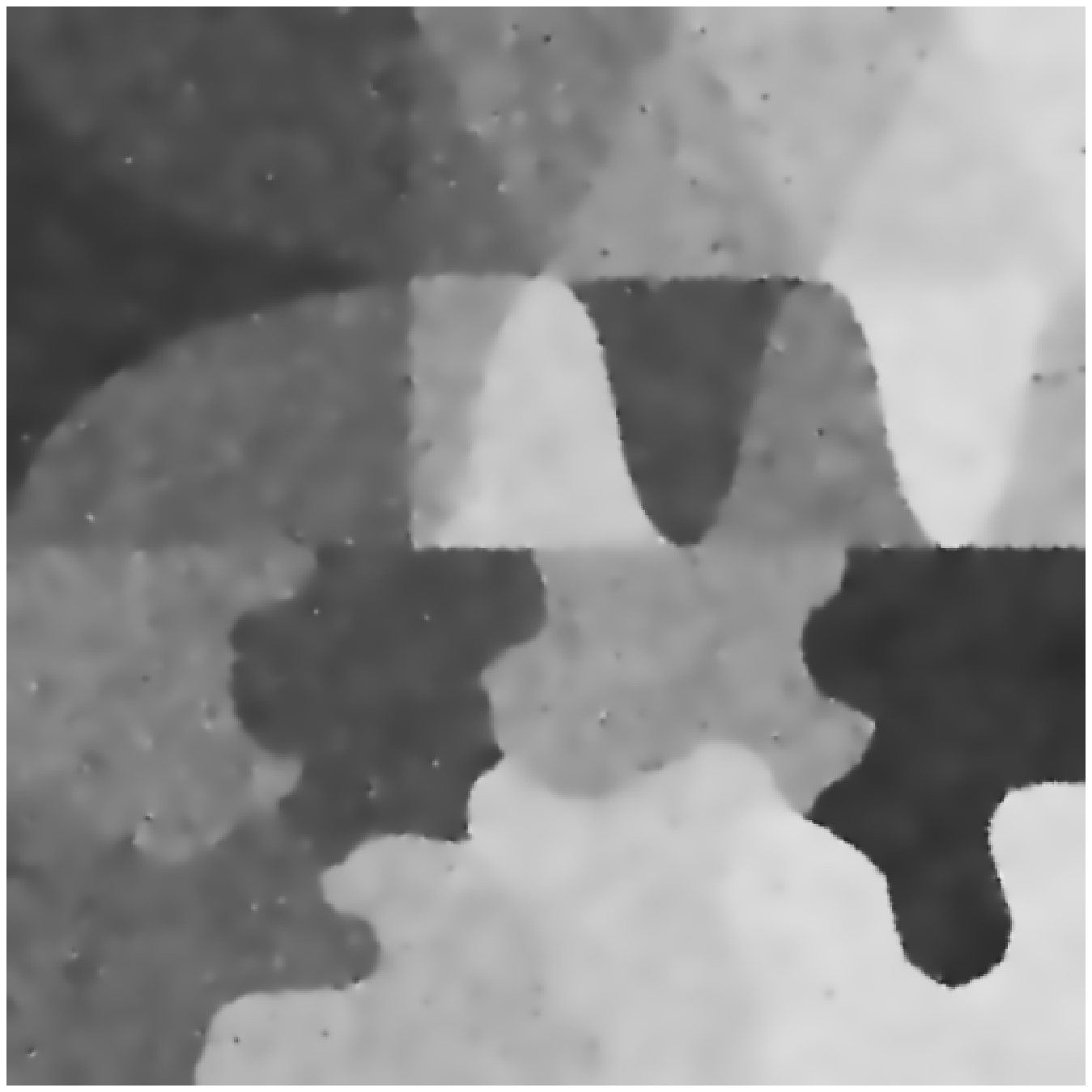}   
                     \caption{}        
                \label{fig:3e}
       \end{subfigure}%
           \begin{subfigure}[b]{0.22\textwidth}           
                \includegraphics[scale=0.19]{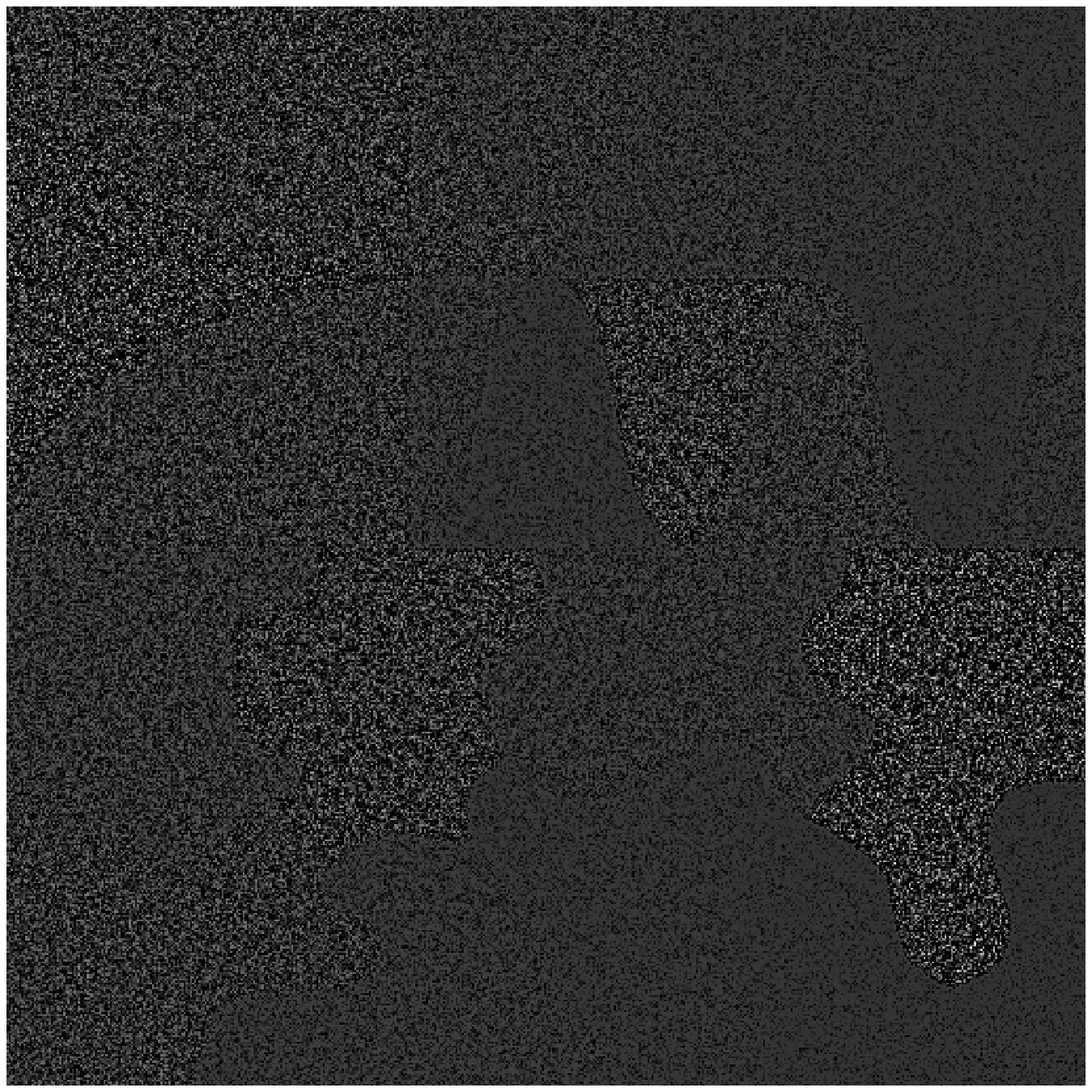}   
                       \caption{}      
                \label{fig:3e_ratio}
       \end{subfigure}%
                  \begin{subfigure}[b]{0.22\textwidth}           
                \includegraphics[scale=0.22]{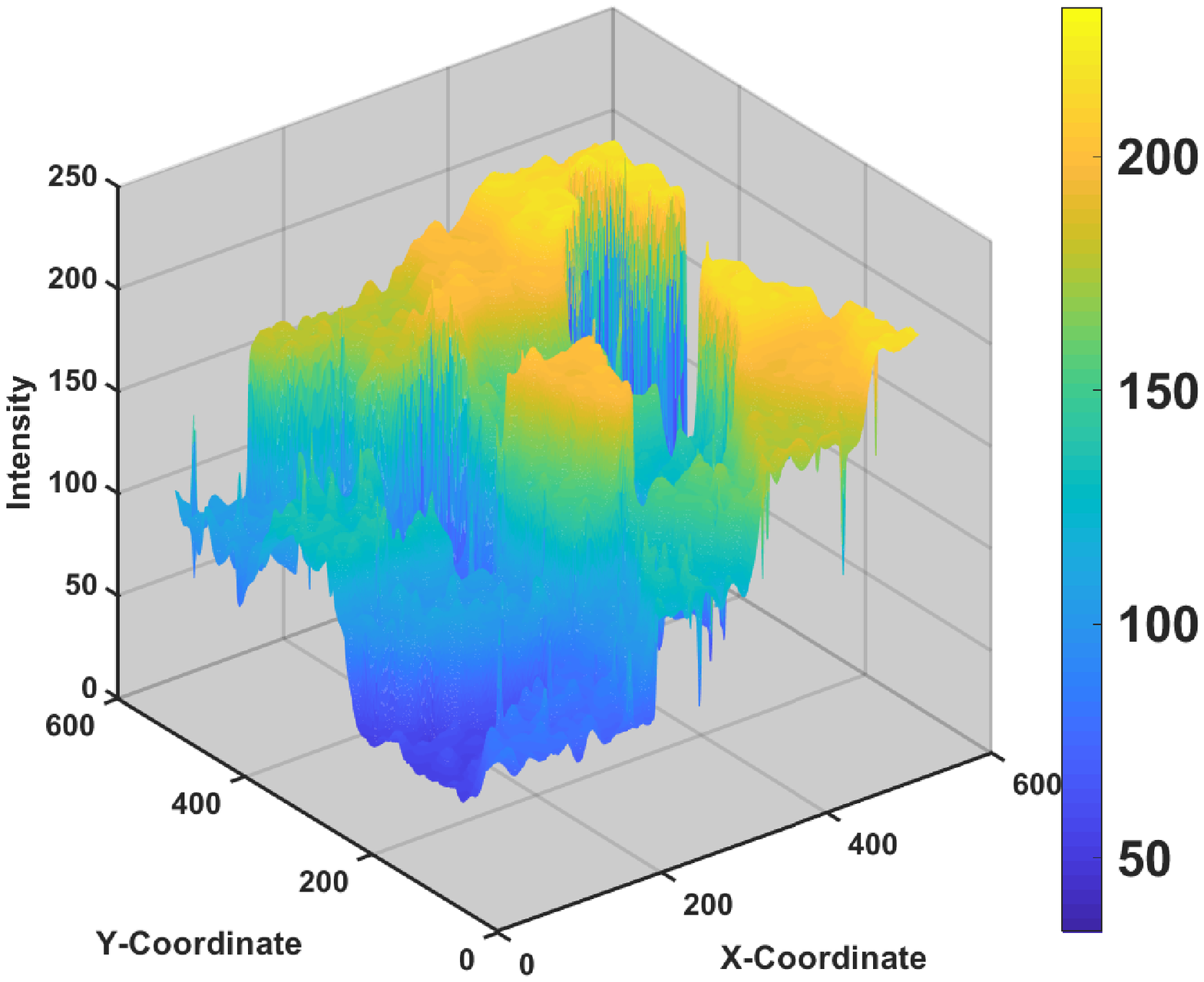}   
                        \caption{}     
                \label{fig:3e_surf}
       \end{subfigure}%

        \begin{subfigure}[b]{0.22\textwidth}           
                \includegraphics[scale=0.19]{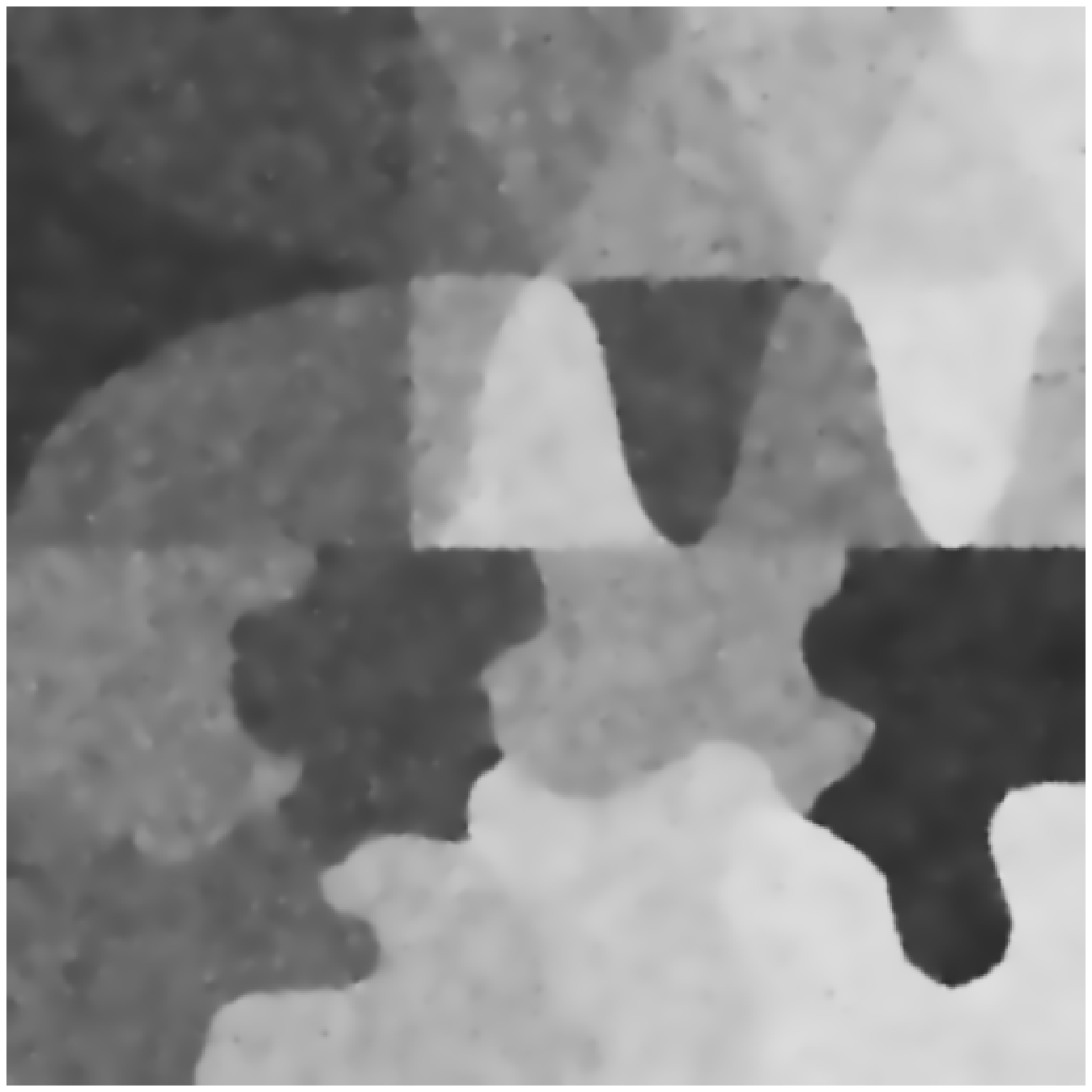}   
                      \caption{}       
                \label{fig:3f}
       \end{subfigure}%
           \begin{subfigure}[b]{0.22\textwidth}           
                \includegraphics[scale=0.19]{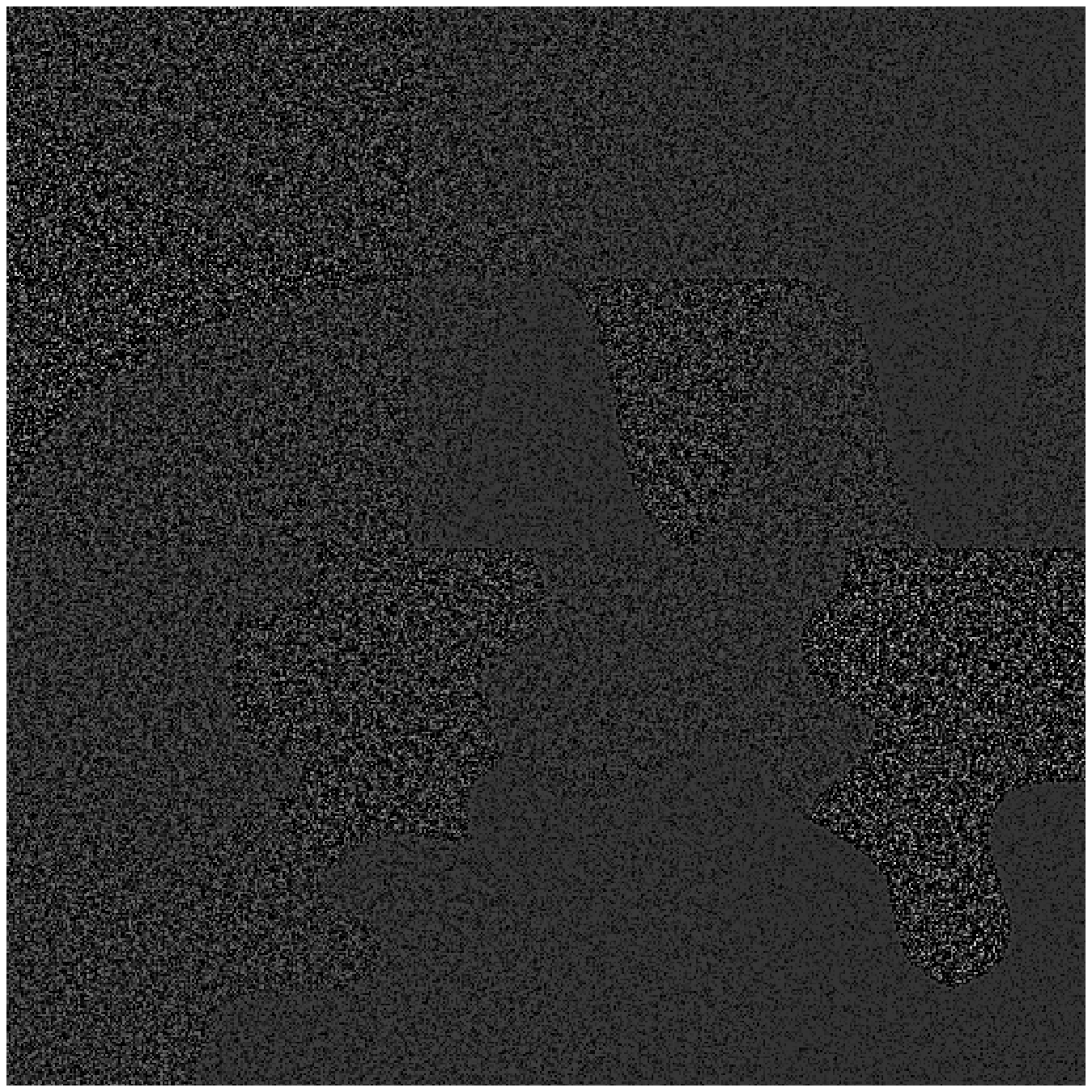}  
                      \caption{}        
                \label{fig:3f_ratio}
       \end{subfigure}%
                  \begin{subfigure}[b]{0.22\textwidth}           
                \includegraphics[scale=0.22]{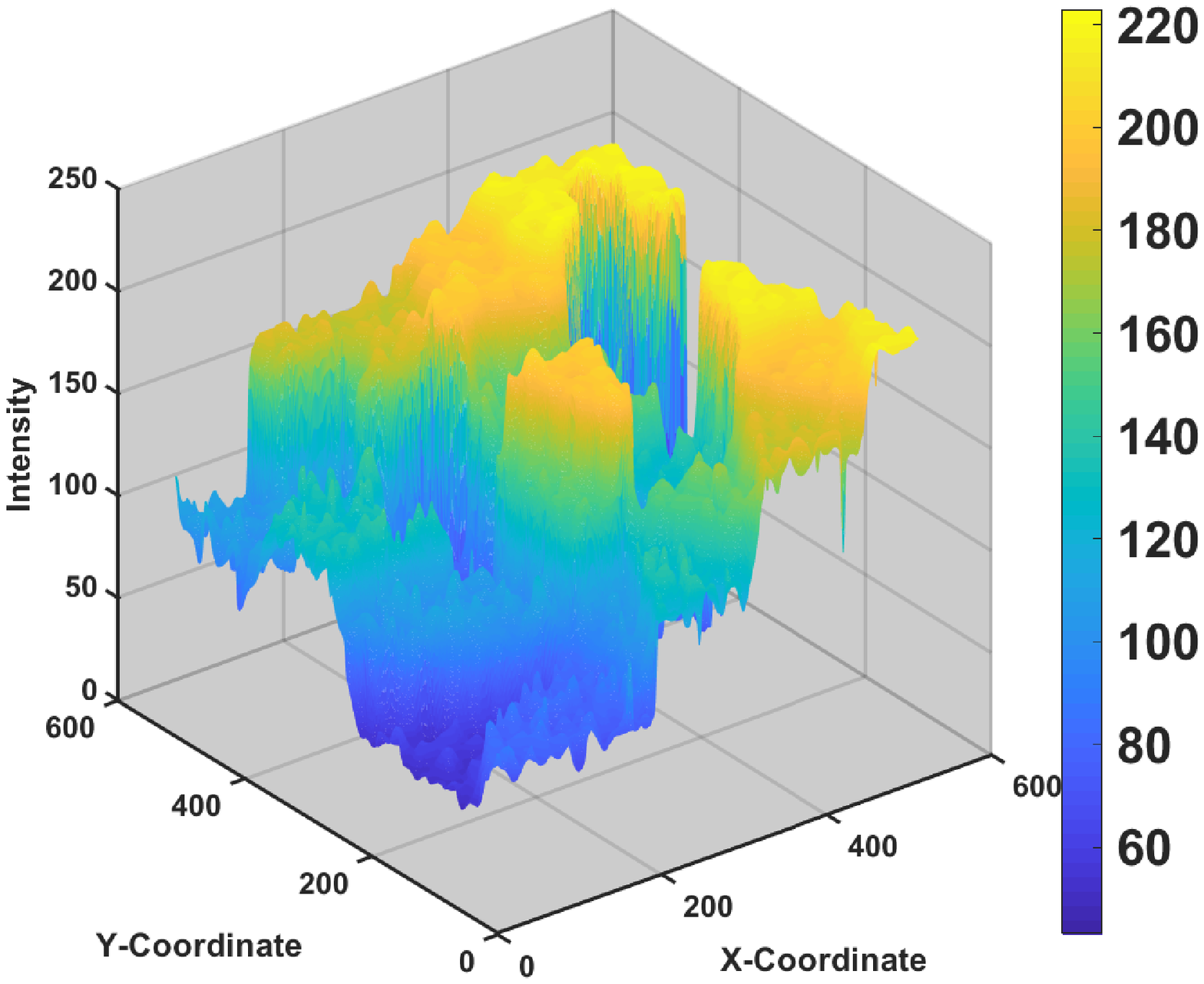}  
                       \caption{}       
                \label{fig:3f_surf}
       \end{subfigure}%
       
 \caption{ A $512 \times 512$  mosaic image corrupted by Gaussian noise with $\sigma=100$ and restored by different models. (a) Original (b) Noisy (c) 3D surface plot of Noisy image (d-f) Cao; $\lambda=20, K=6$ (g-i) SYS; $\lambda=0.1, K=5$ (j-l) ACPDE; $\kappa=1, \nu=1, k=5$ (m-o) Proposed; $\alpha=2, \beta=1, k=5, \nu=1.$ }\label{fig:mosaic_100}
\end{figure}

For further study on the enhancement properties of the proposed model, figure \ref{fig:mosaic_100} depicts the results for Mosaic test images that contaminated by very high-level additive Gaussian noise with $\sigma=100$. Here we present the denoised images, their ratio images, and 3D surface plots of the denoised images. From the first column of the figure, it can be observed that our model works better in terms of noise removal as compared to other models. From the 3D surface plots of the denoised images, one can see that our model leave fewer fluctuations in comparison to the other models. Moreover, from the ratio images, it indicates that the proposed model not only removes noise efficiently but also preserves the fine structure as compared to the other models.

In addition to improved qualitative performance, the quantitative comparisons of different considered methods with different test images and noise levels are shown in Table \ref{tab:table1}, in terms of PSNR and MSSIM values. From the higher values of PSNR and MSSIM, we can observe that the proposed TCPDE model is superior to other existing models consider for the comparison.

Summarizing all the above numerical experiments reveals the better performance of the present model. This model well restored the grayscale images,  as well as preserves the crucial features of the images.


\begin{table}
\centering
\caption{Comparison of MSSIM and PSNR values of denoised images using various approaches for images corrupted by a additive Gaussian noise with different noise level ( standard deviation:$\sigma$)}
\label{tab:table1}
\scalebox{0.85}{
\begin{tabular}{ccrrrrrrrrrr}
\toprule
    \multirow{2}[4]{*}{Image} & \multirow{2}[4]{*}{ $\sigma$ } & \multicolumn{2}{c}{TDM\cite{ratner2007image}} & \multicolumn{2}{c}{Cao \cite{cao2010class}} &\multicolumn{2}{c}{SYS\cite{sun2016class}}& \multicolumn{2}{c}{ACPDE\cite{acpde2019rgate}} & \multicolumn{2}{c}{Proposed} \\
\cmidrule(r){3-4}
		\cmidrule(r){5-6}
		\cmidrule(r){7-8}
		\cmidrule(r){9-10}    
		\cmidrule(r){11-12}  
		      &       & \multicolumn{1}{c}{MSSIM} & \multicolumn{1}{c}{PSNR} & \multicolumn{1}{c}{MSSIM} & \multicolumn{1}{c}{PSNR} & \multicolumn{1}{c}{MSSIM} & \multicolumn{1}{c}{PSNR} & \multicolumn{1}{c}{MSSIM} & \multicolumn{1}{c}{PSNR} & \multicolumn{1}{c}{MSSIM} & \multicolumn{1}{c}{PSNR}\\
    \midrule
Peppers     & 20       & 0.8865 & 29.64 & 0.8878 & 29.89 & 0.8826 & 29.85 & 0.9165 & 31.04 & \textbf{0.9375} & \textbf{31.33} \\
            & 40       & 0.8215 & 27.12 & 0.8273 & 27.27 & 0.8264 & 26.88 & 0.8506 & 27.80 & \textbf{0.8882} & \textbf{28.24} \\
            & 60       & 0.7655 & 25.30 & 0.7718 & 25.41 & 0.7397 & 24.96 & 0.7985 & 25.66 & \textbf{0.8344} & \textbf{26.02} \\
            & 80       & 0.7148 & 23.53 & 0.7225 & 23.61 & 0.6937 & 23.23 & 0.7520 & 23.73 & \textbf{0.7918} & \textbf{24.05} \\
            & 100      & 0.6523 & 21.84 & 0.6600 & 21.89 & 0.6255 & 21.63 & 0.6888 & 21.98 & \textbf{0.7357} & \textbf{22.23} \\
            &          &        &       &        &       &        &       &        &       &                 &                \\
Tree        & 20       & 0.7864 & 26.80 & 0.7802 & 27.24 & 0.7809 & 27.19 & 0.8360 & 27.99 & \textbf{0.8416} & \textbf{28.10} \\
            & 40       & 0.7500 & 24.26 & 0.7046 & 24.41 & 0.7065 & 23.98 & 0.7633 & 24.86 & \textbf{0.7689} & \textbf{24.99} \\
            & 60       & 0.6893 & 22.35 & 0.6490 & 22.42 & 0.6371 & 22.04 & 0.7058 & 22.68 & \textbf{0.7081} & \textbf{22.79} \\
            & 80       & 0.6287 & 20.63 & 0.5919 & 20.71 & 0.5674 & 20.36 & 0.6480 & 20.89 & \textbf{0.6495} & \textbf{21.01} \\
            & 100      & 0.5769 & 19.20 & 0.5451 & 19.29 & 0.5109 & 19.04 & 0.5974 & 19.43 & \textbf{0.5999} & \textbf{19.51} \\
            &          &        &       &        &       &        &       &        &       &                 &                \\           
Mosaic      & 20       & 0.9516 & 33.64 & 0.9164 & 32.54 & 0.9258 & 32.59 & 0.9680 & 34.52 & \textbf{0.9884} & \textbf{34.53} \\
            & 40       & 0.9168 & 29.54 & 0.8666 & 28.81 & 0.8361 & 28.22 & 0.9436 & 29.79 & \textbf{0.9707} & \textbf{30.03} \\
            & 60       & 0.8838 & 26.06 & 0.8267 & 25.63 & 0.7838 & 25.08 & 0.9153 & 26.17 & \textbf{0.9324} & \textbf{26.26} \\
            & 80       & 0.8449 & 23.24 & 0.7818 & 22.98 & 0.6950 & 22.56 & 0.8839 & 23.33 & \textbf{0.9053} & \textbf{23.43} \\
            & 100      & 0.8096 & 21.08 & 0.7455 & 20.91 & 0.6889 & 20.62 & 0.8534 & 21.16 & \textbf{0.8784} & \textbf{21.23} \\
            &          &        &       &        &       &        &       &        &       &                 &                \\            
Brick       & 20       & 0.8622 & 27.28 & 0.8634 & 27.37 & 0.8635 & 27.22 & 0.8610 & 27.31 & \textbf{0.8645} & \textbf{27.42} \\
            & 40       & 0.7335 & 24.66 & 0.7419 & 24.87 & 0.7391 & 24.66 & 0.7326 & 24.97 & \textbf{0.7433} & \textbf{25.00} \\
            & 60       & 0.6505 & 23.30 & 0.6578 & 23.44 & 0.6495 & 23.15 & 0.6638 & 23.64 & \textbf{0.6691} & \textbf{23.72} \\
            & 80       & 0.5839 & 22.08 & 0.5911 & 22.20 & 0.5768 & 21.88 & 0.5992 & 22.31 & \textbf{0.6039} & \textbf{22.35} \\
            & 100      & 0.5244 & 21.04 & 0.5306 & 21.12 & 0.5145 & 20.87 & 0.5347 & 21.15 & \textbf{0.5405} & \textbf{21.22} \\
            &          &        &       &        &       &        &       &                &                 &                \\            
Aerial      & 20       & 0.8981 & 26.39 & 0.8969 & 26.38 & 0.8978 & 26.31 & 0.8999 & 26.48 & \textbf{0.9004} & \textbf{26.49} \\
            & 40       & 0.7789 & 23.85 & 0.7800 & 23.85 & 0.7766 & 23.76 & 0.7820 & 23.87 & \textbf{0.7840} & \textbf{23.90} \\
            & 60       & 0.6776 & 22.56 & 0.6827 & 22.56 & 0.6792 & 22.47 & 0.6857 & 22.57 & \textbf{0.6910} & \textbf{22.60} \\
            & 80       & 0.5893 & 21.56 & 0.5962 & 21.57 & 0.5912 & 21.48 & 0.5968 & 21.57 & \textbf{0.5973} & \textbf{21.57} \\
            & 100      & 0.5283 & 20.79 & 0.5367 & 20.79 & 0.5319 & 20.71 & 0.5372 & 20.79 & \textbf{0.5380} & \textbf{20.80} \\
            &          &        &       &        &       &        &       &                 &                &                \\  
Moon        & 20       & 0.6763 & 29.58 & 0.6676 & 29.49 & 0.6790 & 29.69 & 0.6815 & 29.75 & \textbf{0.6817} & \textbf{29.80} \\
            & 40       & 0.5903 & 27.71 & 0.5860 & 27.56 & 0.5882 & 27.54 & 0.5960 & 27.90 & \textbf{0.5973} & \textbf{27.93} \\
            & 60       & 0.5606 & 26.71 & 0.5564 & 26.57 & 0.5519 & 26.38 & 0.5615 & 26.84 & \textbf{0.5663} & \textbf{26.97} \\
            & 80       & 0.5354 & 25.70 & 0.5318 & 25.59 & 0.5279 & 25.40 & 0.5357 & 25.83 & \textbf{0.5416} & \textbf{25.97} \\
            & 100      & 0.5206 & 24.69 & 0.5183 & 24.61 & 0.5160 & 24.49 & 0.5209 & 24.77 & \textbf{0.5241} & \textbf{24.89} \\
\bottomrule
\end{tabular}
}
\end{table}
\clearpage
\section{Conclusion}
\label{sec:conclusion}
A non-linear coupled telegraph diffusion equation based model for image denoising is introduced in this paper. Here we have taken telegraph equation for the image variable as well as for the edge variable, which improves the present model over the existing coupled PDE based models. Mathematical analysis of the model has been carried out using Banach's fixed point theorem. Also, to validate the effectiveness of the proposed model, numerical experiments are carried out using different standard test images.  Qualitative and quantitative studies, in terms of MSSIM, PSNR, and visual quality, confirm that the proposed model exhibits better performance than the single telegraph-diffusion based models as well as coupled partial differential equation models. This system not only removes the noise but also reduces the staircase artifacts and improves the performance of the filtering, even in low SNR images. Overall, the proposed model is an important addition in the field of diffusion equation based image restoration.
\thispagestyle{empty}

\bibliographystyle{unsrt}

\end{document}